\documentclass[11pt]{amsart}
\newcommand{\sekshun}[1]                
        {                               
        \section{#1}                    
        \markboth{#1 \hfill}{#1 \hfill} 
        }                               
\usepackage{amsfonts, amssymb, curves, eucal,pxfonts, txfonts}
\usepackage{eepic}
\usepackage{pstricks,pst-node, pst-tree}
\usepackage{ifthen}
\newcounter{mycount}

\renewenvironment{proof}{\hfill\newline  Proof:}{\qed\newline\newline}
\newcommand{\R}{\mathbb{R}}
\newcommand{\N}{\mathbb{N}}

\newcommand{\Z}{\mathbb{Z}}
\newcommand{\cT}{\mathbb{T}}
\newcommand{\E}{\mathbb{E}}
\newcommand{\T}{\mathcal{T}}

\newcommand{\A}{\mathcal{A}}

\newcommand{\C}{\mathcal{C}}
\renewcommand{\L}{\mathcal{L}}

\newcommand{\F}{\mathcal{F}}

\newcommand{\n}{{\widetilde n}}
\newcommand{\phiphi}{\Phi}
\renewcommand{\phi}{\varphiup}
\renewcommand{\psi}{\psiup}
\renewcommand{\gamma}{\gammaup}
\renewcommand{\alpha}{\alphaup}
\renewcommand{\beta}{\betaup}
\renewcommand{\tau}{\tauup}
\renewcommand{\sigma}{\sigmaup}
\renewcommand{\rho}{\rhoup}
\renewcommand{\omega}{\omegaup}
\renewcommand{\lambda}{\lambdaup}

\newcommand{\tT}{\tilde{T}}

\newcommand{\ds}{\displaystyle} 
\newcommand{\larr}{\left( \begin{array}{c}}
\newcommand{\rarr}{\end{array} \right) }
\newcommand{\fiveands}{&&&&&}
\newcommand{\tenands}{&&&&&&&&&&}
\newcommand{\allands}{&&&&&&&&&&&&&&&&&&&&&&&&&&&&&&&&&&&&&}
\newcommand{\os}{\overline\sigma}
\newcommand{\oa}{\overline\alpha}
\newcommand{\ob}{\overline\beta}

\newcommand{\arrow}{\rightarrow}

\newcommand{\inv}{\varprojlim}

\begin{document}

\newtheorem{theorem}{Theorem}

\newtheorem{corollary}[theorem]{Corollary}
\newtheorem{lemma}[theorem]{Lemma}
\newtheorem{prop}[theorem]{Proposition}
\newtheorem{example}{Example}

\title{The branch locus for one-dimensional Pisot tiling spaces}
\author{Marcy Barge, Beverly Diamond and Richard Swanson}

 \maketitle
 \begin{abstract} If $\varphi$ is a Pisot substitution of degree $d$, then the
inflation and substitution homeomorphism $\phiphi$ on the tiling space $\T_\phiphi$
factors via geometric realization onto a
d-dimensional  solenoid.  Under this realization,  the collection of
$\phiphi$-periodic asymptotic tilings corresponds  to a finite set that  projects 
onto the \underline{ branch locus} in a
d-torus.  We prove that if two such tiling spaces are homeomorphic, then the
resulting branch loci are the same up to the action of certain affine maps on
the torus.     
\end{abstract} 

\sekshun{Introduction}
 In this paper we
introduce the branch locus, a new topological invariant for one-dimensional Pisot substitution
tiling spaces.  A  substitution  on 
 $n$ letters  is a map from an alphabet $\A= \{1, 2, \ldots, n\}$  into $ \A^*$,
where $\A^*$ is the collection of finite and nonempty words from $\A$. The  \emph{
abelianization matrix }  of $\phi$ is defined as $A = (a_{ij})$, where $a_{ij} =
$ number of occurrences of
$i$ in $\phi(j)$.  The substitution   $\phi$ is  \emph{ Pisot } provided the
Perron-Frobenius eigenvalue $\lambda$ of  $A$  is a Pisot-Vijayaraghavan number
($\lambda > 1$ and all algebraic conjugates of $\lambda$ are strictly inside the
unit circle).  The \emph{ degree } of $\phi$ is the degree of the minimal
polynomial of
$\lambda$.  
 
 Associated with any substitution $\phi$, there is a tiling space
$\T_\phi$ consisting of certain  tilings of the real line;  if
$(\lambda_1, \ldots, \lambda_n)$ is a positive left eigenvector of
$A$, then $\T_\phi$  consists of all tilings of $\R$ by translates of the
prototiles $P_i = [0, \lambda_i]$, $i = 1, \ldots, n$, with the property that
the word spelled out by any finite patch of consecutive tiles in the tiling
occurs as a factor of $\phi^m(i)$ for some $i \in
\A$ and $m \in \N$.   The substitution $\phi$ is \emph{ primitive} if for some $m
\in \N$, every entry of $A^m$ is strictly positive, and \emph{ aperiodic} if no
element of $\T_\phi$ is periodic under translation.  For the remainder of the paper, all substitutions will be
assumed primitive and aperiodic.

 Define the topology of $\T_\phi$ by stipulating that two tilings are close provided a small translate of
one is identical to the other in  a large neighborhood of the origin.  Under the
assumption that $\phi$ is primitive and aperiodic, $\T_\phi$ is a continuum (a
compact, connected metric space).  For such a $\phi$, \emph{ inflation
and substitution} is the  homeomorphism $ {\Phi}: \T_\phi \arrow \T_\phi$ that replaces each
tile $t + P_i = [t, \lambda_i + t]$ of a tiling $T$ in $\T_\phi$ by the patch
$[\lambda t, \lambda \lambda_i + \lambda t]$ tiled by translates of prototiles
following the pattern of the word $\phi(i)$.  There is also a minimal and
uniquely ergodic $\R$-action on $\T_\phi$, called the \emph{ translation flow},
given by $T = \{T_i\}_{i = -
\infty}^\infty\mapsto T - t := \{ T_i-t\}_{i = - \infty}^\infty$, for
$t \in \R$.

The topology of a substitution tiling space is of interest for a number of
reasons.  Physics provides one source of motivation.  Suppose that $T$ is a
tiling in the substitution tiling space
$\T_\phi$.  Placing an atom at the end of each tile creates a one-dimensional
material which is called a \emph{ quasi-crystal} if its diffraction
spectrum is pure point (the atoms must be `weighted' according to the tile types
they lie in).  Bombieri and Taylor (\cite{bomtay}) proved that if $\phi$ is
Pisot, such a material has a nontrivial discrete component in its spectrum. 
Whether  such a material has pure point spectrum when  $\phi$ is irreducible
unimodular Pisot (i.e., $\phi$ is Pisot with ${\rm degree}(\lambda)= d = n = |\A|$, and
$\det(A) =  \pm 1$) remains an open question.  Lee, Moody and Solomyak
(\cite{moodysol}) have proved that the diffraction spectrum of the material is
pure point  if and only if the dynamical spectrum of the translation flow on
$\T_\phi$ is  pure discrete.  It follows from \cite{bargeswanson2}  that if $\T_\phi$
and $\T_\psi$ are homeomorphic tiling spaces, then the tiling flow on $\T_\phi$
is pure discrete if and only if the tiling flow on $\T_\psi$ is pure discrete. 
That is, the question of whether or not a one-dimensional material built from a
substitution has pure point diffraction spectrum is a \emph{ topological} question
about the corresponding tiling space.

Substitution tiling spaces also arise in the study of hyperbolic attractors. 
R.F. Williams (\cite{williams}) proved that every hyperbolic one-dimensional
attractor is topologically conjugate with the shift map on the inverse limit of
an expanding endomorphism of a branched one-manifold, and, with minor
restrictions on the map of the branched one-manifold, all such inverse limits
can be realized as hyperbolic attractors. More recently, Anderson and Putnam
(\cite{ap}) proved that inflation and substitution on a one-dimensional
substitution tiling space is conjugate with the shift on the inverse limit of an
expanding endomorphism of a branched one-manifold.  As a consequence, every
orientable hyperbolic one-dimensional attractor is either a substitution tiling
space,  for which the underlying manifold is branched,  or a classical solenoid,
for which the underlying manifold is the circle.  Modeling an attractor as a
tiling space provides a much clearer view of its global topology than one gets
from considering an inverse limit description:  moving along an arc component in
the attractor is simply translating a tiling, and the patterns of consecutive
tiles determine the recurrence properties of the translates.

Although the ``inverse limit on branched manifolds'' description of tiling spaces
will not play an explicit role in this paper, the intuitive content  of our main
result, and the rationale for the terminology ``branch locus'' that we introduce,
has its origin in that description.  There are, of course, no actual branch
points in a tiling space:  in the one-dimensional case, every point has a
neighborhood that is homeomorphic with the product of an arc and a Cantor set. 
Nevertheless, an inverse limit description of the tiling space gives a sequence
of approximating branched one-manifolds.  In the limit, the ghost of the
branches can be observed in the existence of asymptotic composants: two
distinct tilings $T, T'$ are \emph{ asymptotic} provided $d(T - t, T' - t) \arrow
0$ as $t \arrow \infty$ or as $t \arrow -\infty$.  The arc components of
asymptotic tilings are called \emph{ asymptotic composants}.  Partially sewing up such
asymptotic composants results in a space that does have branching:  the new,
branched, space corresponds to an inverse limit on a branched manifold with 
periodic branch points.  

Although this ``ghost branching" in $\T_\phi$ seems to have no clear location, we
will see that, in case the substitution is Pisot, the branching occurs in
well defined relative geometrical patterns.  The appropriate underlying geometry
is that of the $d$-dimensional torus, where $d$ is the degree of $\lambda$.  Our
main result is that if
$\T_\phi$ and $\T_\psi$ are homeomorphic tiling spaces, then their branch loci,
nonempty finite sets of points in the $d$-torus that we define in \S\ref{geom},
are equal modulo the action of a certain collection of affine endomorphisms of
the torus.  Thus the branch locus becomes a topological invariant.  We
illustrate this by distinguishing pairs of tiling spaces that are otherwise
difficult to tell apart (see \S\ref{examples}). 

The idea for considering ``branching'' in one-dimensional tiling spaces arose in 
discussions the first author had with S{\o}ren Eilers regarding the topological content of the 
$K_0$-group of the Matsumoto algebra associated with a substitution. 
Eilers, Restorff and Ruiz (\cite{eilresruiz}) have shown that this (ordered) group is also a complete invariant of the Matsumoto algebra and, consequently, the Matsumoto algebra is a topological invariant of the tiling space.


How is the Matsumoto algebra (equivalently, its $K_0$-group) reflected in the topology of the 
tiling space? In \S\!\!~\ref{pisotpartsection} we show that the branch locus provides a partial answer. We use the 
branch locus to define the ``Pisot part of the augmented dimension group'', an ordered group that is a flow equivalence invariant of the substitution and an ordered subgroup of the augmented cohomology group of the tiling space, which,
in turn, is closely related to the Matsumoto $K_0$-group (see \cite{carlseneilers1}, \cite{carlseneilers2} for a description of the Matsumoto $K_0$-group in the substitutive setting and \cite{bargesmith} for an account of the relationship between the Matsumoto $K_0$-group of a substitutive system and the augmented cohomology of the associated tiling space). 

\sekshun{Geometric Realization} \label{geom}  For convenience, we will  use the 
``strand space'' model for the tiling space (see
\cite{marcyjarek} and \cite{bakermarcyjarek}) which we recall now.  Let $\phi$ be a 
Pisot substitution of degree $d$ on $n$ letters with abelianization $A$
and Perron-Frobenius eigenvalue $\lambda$. 
There is an (unique) $A$-invariant decomposition $\R^n = V \oplus W$ such that 
$V$ contains a right Perron-Frobenius eigenvector $\omega$ associated with  
$\lambda$ and $\dim V=d$\footnote{$V$ is the kernel of $m(A)$, if   $m(x)$
denotes the minimal polynomial of $\lambda$. }. 
We can always choose rational bases for $V$ and $W$. There is a further  
$A|_V$-invariant splitting $V=\E^u\oplus\E^s$ obtained by letting $\E^s$ be the space orthogonal to a left Perron-Frobenius eigenvector of $A|_V$ and $\E^u$ be the span of $\omega$.  
Let $\mathrm{pr}_V:\R^n\to V,\ \mathrm{pr}_s:V\to\E^s,\ \mathrm{and}\ \mathrm{pr}_u:V\to\E^u$ denote the projections, resp., along $W$, $\E^u$ and $\E^s$, and let $\Gamma$ denote the $A$-invariant lattice  
$\mathrm{pr}_V \Z^n$.

If  $e_i, i = 1, \ldots, n$,  are the
standard basis vectors in $\R^n$, define  $\varv_i := \mathrm{pr}_V(e_i)$, and let $\sigma_i
:= \{t\varv_i: 0 \leq t \leq 1\} $ denote the oriented segment representing $\varv_i$. 
Even if $\sigma_i = \sigma_j$ for some $i
\neq j$, we wish to distinguish between these segments: we call 
$\sigma_i$  a \emph{ (labeled) segment of type $i$}.  An oriented broken line
$\gamma = \{\sigma_{i_k} + x_k\}$, $x_k \in V$, consisting of a collection of
translated copies of the basic segments meeting tip-to-tail and with connected
union, will be called a  \emph{ strand}.  We denote the  space of bi-infinite
strands in $V$ by
$$\F:= \{\gamma: \gamma  \textrm{ is a bi-infinite  strand in } V\}.$$  The
substitution $\phi$ induces the \emph{ inflation and substitution map} $ {\Phi}: \F
\arrow \F$ as follows:  for each edge (translated segment) $\sigma_{i_k} + x_k$
in strand $\gamma$, replace that edge by the collection of edges $\sigma_{j_1} +
Ax_k$, 
$\sigma_{j_2} + Ax_k + \varv_{j_1}$, $\ldots$, $\sigma_{j_\ell} + Ax_k + \varv_{j_1} +
\cdots + \varv_{j_{\ell-1}}$, where $\phi(i_k) = j_1 j_2 \cdots j_\ell$.  That is,
$\Phi(\gamma)$ is obtained by applying the linear map
$A$ to each edge of $\gamma$, then breaking up the image into translated
segments following the pattern determined by applying $\phi$ to the type of the
edge.

For $R > 0$, let $\F^R$ denote the subset of $\F$ consisting of those strands
all of whose edges are within distance $R$ of $\E^u$.  There is then an $R_0$ so
that $$\F_\phi := \bigcap_{n \geq 0} {\Phi}^n(\F^R)$$ is independent of $R
\geq R_0$.  The set $\F_\phi$ has a natural metric topology in which $\gamma$ and
$\gamma'$ are close if a small translate of $\gamma$ by a vector in $V$ lines
up exactly with $\gamma'$, segment types being considered, in a large
neighborhood of the origin.  With respect to this topology, 
 $ {\Phi}$ and the \emph{ translation flow}, 
\[
\gamma = \left\{\sigma_{i_k} + x_k\right\} \mapsto
\gamma - t := \left\{\sigma_{i_k} + x_k - t \omega\right\},
\] 
are continuous.

In some cases there may be a few translation orbits in $\F_\phi$ that correspond
to strands with acceptable heads and tails but that are joined in an unnatural
way (this happens, for instance, when there are letters $a$ and $b$ with $\phi(a)
= a \cdots$, $\phi(b) = \cdots b$, but  the word $ba$ never occurs in $\phi^m(i)$,
$i \in \A$, $m \in
\N$).  We eliminate these chimeras by defining $\T_\phi$ to be the
$\omega$-limit set, under  translation flow, of any $\gamma \in
\F_\phi$.  That is, for any $\gamma \in \F_\phi$, $$\T_\phi := \bigcap_{T > 0}
\textrm{cl} \{\gamma - t: t \in [T,
\infty)\}.$$

 Inflation and substitution, $ {\Phi}:
\T_\phi \arrow \T_\phi$,  is a homeomorphism from $\T_\phi$
 onto itself,  and the translation flow  on $\T_\phi$ is minimal and uniquely
ergodic.
 
One reason for using the rather elaborate strand space model of the tiling space
is that it permits a simple and concrete definition  of geometric realization. 
One could  factor the tiling dynamics onto those of a solenoid, by choosing the
lattice $\Gamma$ in $V$ and mapping $\gamma  \in \T_\phi$ to $$( \min I_0 +
\Gamma, \min I_1 + \Gamma,
\ldots) \in \inv F_{A|_V}, $$  where $F_{A|_V} : V/\Gamma \arrow V/\Gamma$ is
defined by $F_{A|_V}(v + \Gamma) = Av + \Gamma$, $I_k$ is any edge in
$ {\Phi}^{-1}(\gamma)$ and $\min I_k$ is its initial vertex.  This does give a well
defined surjection of $\T_\phi$ onto the
$d$-dimensional solenoid $\inv F_{A|_V}$ that semiconjugates inflation
and substitution with the shift, as well as    translation flow with a Kronecker
action.  However,   to  maximize the size of the factor, we will define a
coarser, more natural lattice.  Toward this end, denote the collection of \emph{
return vectors} by \begin{eqnarray*} \Theta(i) := & \{v \in
\Gamma: \textrm{ there exists } \gamma \in \T_\phi \textrm{ containing edges } I,
I', \\ &\hspace*{1.5in}
\textrm{ each of type } i, \textrm{ with } I' = I + v\}.\end{eqnarray*}

{\parindent= 0pt It is not} difficult to show that the subgroup
$$\Sigma_\infty :=  \left<  \bigcup_{k \in \Z}(A|_V)^k \Theta(i)
 \right>$$ of $V$ generated by $ \bigcup_{k \in \Z}(A|_V)^k \Theta(i)$ is
independent of $i \in \A$.  The \emph{ return lattice} $$\Sigma :=
\Sigma_\infty \cap \Gamma$$ is invariant under $A$ and also of rank
$d$.
 
 Note:  In case $A$ is irreducible ($d=n$), then
$\Sigma = \Gamma = \Z^n$.  In \S\!\!~\ref{examples}, we will see an example where
$\Sigma$ is strictly coarser than $\Gamma$.
 
 To define geometric realization onto a solenoid determined by
$\Sigma$, we must determine some appropriate translations.  We begin by noticing
that for each $i, j \in \A$, there is a well defined ``transition vector''
$\varw_{ij} \in \Gamma/\Sigma$ so that if $\gamma \in
\T_\phi$ and $I$ and $J$ are edges of $\gamma$ of types $i$ and $j$ resp.,
then $$(\min J - \min I)\textrm{ mod } \Sigma = \varw_{ij}.$$  Then, since $\varw_{ij} +
\varw_{jk} + \varw_{ki} = 0$ for all $i, j, k
\in \A$,  there are $u_i \in \Gamma/\Sigma$, for $i \in \A$, so that for all $i,
j$, 
 \[\varw_{ij} = u_j - u_i.\]
 We may then normalize the $u_i$ so that $$A u_i = u_{i'}$$ for all $i
\in \A$, where $i'$ denotes the initial letter of $\phi(i)$.  Letting
$F_{A|_V}$ denote the map induced by $A$ on the torus $V/\Sigma$, we have a
well defined map $\varg_\phi: \T_\phi \arrow \inv F_{A|_V}$ (we will always use the lattice $\Sigma$ in place of $\Gamma$ in the definition of $F_{A|_V}$) given by 
 \[ 
\varg_\phi(\gamma) = (
 \min I_0 + u_{i_{0}}, \min I_1 + u_{i_{1}}, \ldots)  
\] 
where $I_k$ is any edge
of $ {\Phi}^{-k}(\gamma)$ of type
$i_k$.  The map  $\varg_\phi$ is called  \emph{geometric realization}.
 
 Note:  The map $\varg_\phi$ depends on the choice of the $u_i \in
\Gamma/\Sigma$, which are not uniquely defined even after the normalization.  
 
 From results of \cite{marcyjarek} and
\cite{bakermarcyjarek}, the map $\varg_\phi$ is boundedly finite-to-one and almost
everywhere $m$-to-$1$, where $m$ is the \emph{ coincidence rank} of
$\phi$, and the tiling flow on $\T_\phi$ has pure discrete spectrum if and only
if   $m = 1$.  Moreover, $\varg_\phi$ is optimal in the sense that any other factoring 
of $ {\Phi}$ onto a solenoidal shift factors through
$\varg_\phi$.\footnote{ For the translation flow on
$\T_\phi$, the Kronecker flow, $( z_{0}, z_{1}, \ldots)
\mapsto (z_{0}- t
\omega, z_{1}- t \lambda^{-1} \omega, \ldots)$ on $\inv F_{A|V}$ is the maximal
equicontinuous factor. }
 
 Geometric realization expresses the underlying solenoidal nature of Pisot
tiling spaces.  But tiling spaces, unlike solenoids, are not homogeneous.  \emph{ Asymptotic tilings} are one kind of inhomogeneity collapsed out by geometric realization. We say  the tilings $\gamma, \gamma'  \in \T_\phi$
are \emph{ forward (backward, {\rm resp.}) asymptotic} provided $d(\gamma - t,
\gamma' - t) \arrow 0$ as $t \arrow \infty$ ($-\infty$,
resp.).\footnote{There is a weakening of asymptoticity, called regional
proximality, with the property that $\varg_\phi(\gamma) = \varg_\phi(\gamma')$ if and only
if $\gamma$ and $\gamma'$ are regionally proximal, see \cite{aus}.}  If $\gamma$
and $\gamma'$ are asymptotic, where $\gamma \neq \gamma'$, there is a unique
$t$ so that
$\gamma - t$ and $\gamma' - t$ are periodic under inflation and substitution.
Moreover, the set   
$\C_\phi$ consisting of those tilings $\gamma\in\T_\phi$ for which there is a 
tiling $\gamma'\ne \gamma$ with $\gamma$ and $\gamma'$  asymptotic and $ {\Phi}$-periodic  is finite and nonempty (see \cite{bargediamond1}). We will call the elements of $\mathcal C_\phi$ \emph{special tilings}.  

Let $L: (V, \Sigma) \arrow (\R^d, \Z^d)$ be a linear isomorphism, and 
 $F_L: V/\Sigma \arrow \cT^d := \R^d/\Z^d$ the induced  isomorphism. Also, let $M=M_\phi$ denote the integer matrix 
representing the linear isomorphism $L\circ A|_V\circ L^{-1}$ in the standard basis, and let $F_M:\mathbb T^d\hookleftarrow$ denote the corresponding toral endomorphism.

  The \emph{branch locus} of $\T_\phi$ is the set $$Br(\phi) := (F_L\circ \pi_0\circ \varg_\phi)(\C_\phi) \subset
\cT^d,$$ where $\pi_0: \inv  F_{A|V} \arrow V/\Sigma$ is projection onto the
$0^{th}$ coordinate.  Note that $F_M(Br(\phi))=Br(\phi)$. Also note that $Br(\phi)$  depends not only on $\phi$ but also on
the choice of $\{u_i\}$ in the definitions of $\varg_\phi$ and the linear
isomorphism
$L$  made in the construction of $\varg_\phi$.  Our main theorem is that this
dependence is limited.

\begin{theorem} \label{main} Suppose that $\phi$ and $\psi$ are primitive,
aperiodic, Pisot substitutions whose tiling spaces are homeomorphic.  Then there
are $d
\times d$ integer matrices $S, T$,   and  $m_0,\ m_1 \in \N$ so that $M_\phi^{m_0} = ST$, 
$M_\psi^{m_1}=TS$, and 
translations $\tau_0,\ \tau_1$ on $\mathbb T^d$ so that 
\[
Br(\psi) = (\tau_0 \circ F_{T})\left(Br(\phi)\right)\quad \text{and}\quad Br(\phi)=(\tau_1\circ F_S)\left(Br(\psi)\right)
\]
\end{theorem}  
  
  In case $\lambda$ is a \emph{ Pisot unit}, that is,  $\det (A|_V) = \pm 1$;
calculations are simplified.
  
  \begin{corollary}\label{pisotunit}  Suppose that $\phi$ and $\psi$ are primitive and  aperiodic
Pisot substitutions with $\lambda_\phi$ a Pisot unit of degree $d$, and that
$\T_\phi$ is  homeomorphic with $\T_\psi$.  Then there are  a $T \in GL(d, \Z)$ and
a translation $\tau$ on $\mathbb T^d$ so that 
\[
Br(\psi) = \tau \circ F_T(Br(\phi)).
\]
  \end{corollary}
  
  Before proving the theorem, we look at a few examples.  
  
 \sekshun{Examples} \label{examples} {\parindent=0pt    

\begin{example}\end{example}
Consider the substitution $\ds \phi:  \begin{cases} 1 \quad \arrow & 121
\\ 2 \quad \arrow & 312 \\3 \quad \arrow & 213
\end{cases} .\quad $  
In this instance, 

\vspace{.1in}

$\displaystyle A =  \left( \begin{array}{ccc} 2 & 1 & 1
\\ 1 & 1 & 1 \\0 & 1 & 1\end{array} \right)$, $\lambda = 3$ and $V = \E^u =
\left\{ t \left( \begin{array}{c} 3  \\ 2  \\1 
\end{array} \right): t \in \R \right\}$.  We have $\varv_i = \mathrm{pr}_Ve_i = \displaystyle \frac{1}{6} 
\left( \begin{array}{c} 3  \\ 2  \\1
\end{array} \right)$, which we denote by $v$, so that $\Gamma = \Z v$.  Now
consider the strand $\gamma$, fixed under $ {\Phi}$,
\[
\gamma = \left\{ \ldots, \sigma_1 - 2v, \sigma_2 - v, \sigma_1,
\sigma_2 + v, \sigma_1 + 2v, \ldots\right\}.
\]  
That is, $\gamma$ follows the pattern
of the fixed word $\cdots 312 . 121 \cdots$ of $\phi$.  Clearly, if $I$
and $I'$ are two edges of $\gamma$ of type 1, then $\min I' - \min I$ is of the
form $2kv$, where $k \in \Z$ (and $k = 1$ occurs).  Thus $$\Sigma = \left\langle
\bigcup_{k \in \Z} (A|_V)^k
\Theta(1)\right\rangle \cap \Gamma = \left< \bigcup_{k \in \Z} (3^k) (2v)\right> \cap
\Gamma = 2 \Gamma = 2 \Z v.$$
 The transition vectors are (mod $\Sigma$) $\varw_{21} = v$, $ \varw_{31} =v$, and
$\varw_{23} = 0$.  Choosing $u_1 = 0$, $u_2 = v$, $u_3 = v$, we have $\varw_{ij} = u_j -
u_i$ and $Au_i = 3u_i = u_{i'} $ (mod  $\Sigma$), where $\phi(i) = i' \cdots$. 
Finally, we map $V/\Sigma$ onto $\R/\Z =
\cT^1$ by $F_L$, where  $L : V \arrow \R$ is defined by $L(tv) =
\frac{t}{2} v$.  Then geometric realization of $\T_\phi$ onto a 3-adic solenoid
is given by $\varg_\phi (\gamma) = ( \min I_{0} - u_{i_{0}},  \min I_{1} - u_{i_{1}},
\ldots)$,  where
$I_k$ is an edge of $ {\Phi}^{-k}(\gamma)$ of type $i_k$ and $\min I_k$ is taken
mod $\Sigma$.  The branch locus of $\phi$ is $(F_L\circ \pi_0 \circ \varg_\phi)(C_\phi)$.
 
 There are procedures for finding asymptotic tilings described in
\cite{bargediamond1} and \cite{bargediamondholton}, and an applet for this
purpose can be found at \cite{eilers1}.  For this $\phi$, there is a single pair
of backward asymptotic tilings (of period 2 under $ {\Phi}$) and a single pair of
forward asymptotic tilings (fixed by $ {\Phi}$) in $\C_\phi = \{\alpha, \alpha',
\beta,
\beta'\}$, where $$\alpha := \{
\ldots, \sigma_2 - v, \sigma_1, \sigma_2 + v, \sigma_1 + 2v, \ldots\}$$ follows
the pattern $\cdots 312. 121 \cdots $,  $$\alpha' := \{
\ldots, \sigma_3 - v, \sigma_1, \sigma_2 + v, \sigma_1 + 2v, \ldots\}$$ follows
the  pattern $\cdots 213. 121 \cdots ,$  $$\beta := \{
\ldots, \sigma_1 - v, \sigma_3, \sigma_1 + v, \sigma_2 + 2v, \ldots\}$$
follows the pattern     $\cdots 121.312 \cdots$ and $$\beta' :=  \{
\ldots, \sigma_1 - v, \sigma_2, \sigma_1 + v, \sigma_3 + 2v, \ldots\}$$ follows
the pattern $\cdots 121.213 \cdots$.  Thus $Br(\phi) = \{0 +
\Z, 1/2 + \Z\}$.  The periodic forward asymptotic tilings are ``halfway around''
the tiling space from the periodic backward asymptotic tilings. 
\qed

Note:  One can show that $\Gamma/\Sigma \simeq \Z/h\Z$, where $h$ is the {\it height}
of the substitution $\phi$ (see \cite{dekking}), if $\phi$ is
constant length (as in Example 1) or, more generally,  $\phi$ is Pisot of
degree 1.  In this sense, the group $\Gamma/\Sigma$ is an extension of the
notion of the height of a constant length substitution to the higher degree
setting.

\begin{example} \label{new}\end{example} We will show that two particular
tiling spaces are not homeomorphic by considering their branch loci.  The
Matsumoto $K_0$-group (\cite{matsumodo}) and the crossing group 
(\cite{bargeswanson1} do not distinguish these spaces.  Also, we have not been
able to distinguish these spaces by  the generalized Bowen-Franks trace
(\cite{bargeswanson1}) or   the  proximality structures described in
\cite{bargediamond2}.  The (unordered) augmented cohomology groups of these spaces 
are isomorphic, but one can show that the Pisot parts of the augmented dimension groups
(see \S\!\!~\ref{pisotpartsection}) are not order isomorphic.

\vspace{.1in}

Consider  the substitutions 
\[
\phi: 
 \begin{cases} a & \to\, aaa^{12} b^{16} a a
\\ b & \to\,  aba^{14} b^{14} ba \end{cases}\quad\text{and}\quad\psi:
\begin{cases} a & \to\, aa aa^{11} b^{16} a a
\\ b & \to\,  aaba^{13} b^{14} ba \end{cases}. 
\]  

\vspace{.1in}

Here $A_\phi = A_\psi = 
\begin{pmatrix} 16 & 16\\ 16 & 16
\end{pmatrix} $, $\varv_\phi = \varv_\psi = \E^u = \left\{t  \left(
\begin{array}{c} 1\\ 1 \end{array} \right): t \in \R\right\}$, $\varv_1 = \varv_2 = v = 
\begin{pmatrix} \frac{1}{2}\\ \frac{1}{2} \end{pmatrix} $ for both
$\phi$ and $\psi$, and $\Sigma_\phi = \Gamma_\phi = \Gamma_\psi =
\Sigma_\psi = \Z v$.  In geometric realization all the $u_i$ can be taken to be
0.  Choosing  $L: tv \mapsto t$,  one gets $Br(\phi) = \{ 1/31 + \Z, 30/31 +
\Z\}$, $Br(\psi) = \{ 2/31 + \Z, 30/31 + \Z\}$, and
$M = M_\phi = M_\psi = (32)$.  So, if $\T_\phi$ is homeomorphic   with
$\T_\psi$, there must be $S = (s)$ and $T = (t)$, where $s, t \in
\Z$,  an $m \in \N$ so that $M^m = ST$, and a translation $\tau$ on $\cT^1 =
\R/\Z$ so that $\tau\circ F_{T}Br(\phi) = Br(\psi)$.  The only possibility for $t$ is
$\pm 2^k$, where $k \in \N$.  Then $\tau\left((\pm 2^k) \{1/31, 30/31\}\right) = \{2/31,
30/31\} (\textrm{mod} 1) $ implies that $(\pm 2^k)(29/31) = (28/31)(\textrm{mod }
1)$, hence that 
$(\pm 2^k) (29) = 28 (\textrm{mod } 31)$, or $\pm 2^{k+1} = 3
(\textrm{mod } 31)$, which is not the case for any $k \in \N$.  It follows from
Theorem \ref{main}  that $\T_\phi$ and $\T_\psi$ are not homeomorphic.
\qed 

\

\begin{example}\end{example} Define  $\phi$:  $\begin{array}{cl} 1 &
\arrow 1122111\\ 2 & \arrow  1211111221 \end{array} $ and $\psi$:
$\begin{array}{cl} 1 & \arrow 1121211
\\ 2 & \arrow  1211112211 \end{array}. $ 

The matrix  $A= A_\phi = A_\psi =  \begin{pmatrix} 5 & 7\\ 2 & 3
\end{pmatrix} $, so $\lambda=\lambda_\phi=\lambda_\psi$ is a Pisot unit and Corollary~\ref{pisotunit} applies.  
As  $\lambda $ has degree 2,  $\Sigma_\phi =
\Sigma_\psi = \Z^2$.  Taking $L $ to be the identity map,  geometrical
realization for both $\T_\phi $ and $\T_\psi$ takes the form $\gamma \mapsto \min
I\ (\textrm{mod } \Z^2)$, where this last $I$ is an edge of $\gamma$.    Each of
$\phi$ and $\psi$ has a pair of backward and a pair of forward asymptotic
tilings,  fixed under inflation and substitution.  Consider, for example, the
tilings for $\phi$ following the patterns 
\[
\cdots 11221111 \dot{1} 22111\cdots
\]
and  $$\cdots 11221111\dot{2}11111221\cdots, $$ where the dot in each word above
indicates the position of the origin.  To locate the corresponding strands (call
them $\gamma_1$ and $\gamma_2$), we seek
$x \in \R^2$ so that $$ {\Phi}(\sigma_1 + x) = \{\sigma_1 + x -e_1,
\sigma_1 + x, \sigma_2 + x + e_1, \ldots\}$$ or, equivalently,
$$ {\Phi}(\sigma_2 + x) = \{\sigma_1 + x -e_1, \sigma_2 + x, \sigma_1 + x + e_2,
\ldots\}.$$  Since $$ {\Phi}(\sigma_1 + x) =  \{\sigma_1 + Ax,
\sigma_1 + Ax + e_1, \ldots\}$$ we must have $Ax = x-e_1$.  This yields $x =
\left( \!\!\begin{array}{r} \frac{1}{3}\\ -\frac{1}{3} \end{array} \right)$, a vertex of both
$\gamma_1$ and $\gamma_2$.  Thus 
\[
\varg_\phi(\gamma_1) = \varg_\phi(\gamma_2) = 
\begin{pmatrix} \,\,\frac{1}{3}\\ \!-\frac{1}{3} \end{pmatrix} =  \begin{pmatrix} \frac{1}{3}\\ \frac{2}{3} \end{pmatrix}
(\textrm{mod } \Z^2) \in Br(\phi).
\] 
Similar calculations yield:  
\[
Br(\phi) =
\left\{
\left( \begin{array}{c} \frac{1}{3}\\ \frac{2}{3} \end{array} \right) + \Z^2,
\begin{pmatrix} \frac{2}{3}\\ \frac{1}{3} \end{pmatrix}  + \Z^2 \right\}
\] 
and 
\[
Br(\psi) = \left\{ \begin{pmatrix} \frac{1}{3}\\ \frac{2}{3} \end{pmatrix} +
\Z^2,  \begin{pmatrix} \frac{1}{2}\\ 0 \end{pmatrix}  + \Z^2\right\}.
\] 
Were $\T_\phi$ and $\T_\psi$ homeomorphic, there would be a 
$T \in GL(2, \Z)$ and a translation
$\tau$ on $\mathbb T^2$ so that $(\tau \circ F_T)(Br(\phi)) = Br(\psi)$. In other words, 

 \begin{eqnarray*} 
 F_T\left[ 
\begin{pmatrix}\frac{1}{3}\\ \frac{2}{3} \end{pmatrix}  
           - \begin{pmatrix} \frac{2}{3}\\ \frac{1}{3} \end{pmatrix}+ \Z^2
\right]
=    \pm \left[ 
 \begin{pmatrix} \frac{1}{3}\\ \frac{2}{3}
 \end{pmatrix}  - \begin{pmatrix} \frac{1}{2}\\ 0
 \end{pmatrix}  + \Z^2 \right].
 \end{eqnarray*} 
That is, $T
\in GL(2, \Z)$ would satisfy $$ \frac{1}{3}\, T\begin{pmatrix} 1\\ 2\end{pmatrix}
- \frac{1}{6}\, \begin{pmatrix} 1\\ 2\end{pmatrix}
\in \Z^2.$$  Since this is not possible, Corollary~\ref{pisotunit} implies that
$\T_\phi$ and $\T_\psi$ are not homeomorphic. \qed

\begin{example}\end{example}

\begin{equation*}
\phi:=\begin{cases}a\to aaa^2b^4cac\\ b\to aba^7b^6c^5bbc\\c\to aaba^6b^6c^3cbc
\end{cases}\qquad
 \psi:= 
\begin{cases}a\to aaab^4caac\\ b\to aaba^5b^5c^5bbbac\\c\to aabbaa^4b^4c^3cbbac
\end{cases}
\end{equation*}
with abelianization
$\ds A=\begin{pmatrix}5&8&8\\4&9&8\\2&6&5\end{pmatrix} $ . The projection along $\begin{pmatrix}2\\0\\ \!\!-1\end{pmatrix}$ onto the Pisot subspace $V$ is given by   
$\ds \mathrm{pr}_V=\begin{pmatrix}0&1&0\\0&1&0\\\frac{1}{2}&\!\!-\!\frac{1}{2}&1\end{pmatrix}$ .
In this case the lattices $\Gamma$ and $\Sigma$ are equal, since return vectors span $\Gamma$ over $\Z$ (the repeated letters $aa,\ bb,\ {\rm and }\ 
cc$  yield return vectors $\varv_i=\mathrm{pr}_V(e_i),\ i=1,2,3$).

\vspace{.1in}
The 
set $\ds \{\varv_1,\varv_1+\varv_2\}=\left\{\begin{pmatrix}0\\0\\\frac{1}{2}\end{pmatrix},\linebreak \begin{pmatrix}1\\1\\0\end{pmatrix}\right\}$ is a basis for $\Sigma$ over $\Z$.
  In these coordinates, $A|_V$ takes the form $\ds M=\begin{pmatrix}5&16\\4&13\end{pmatrix}$.  

\vspace{.1in}

There are one pair of fixed backward asymptotic tilings and two pairs of fixed forward asymptotic tilings for both $\phi$ and $\psi$. We solve the equations for the branch locus points in a manner similar to the foregoing examples: 
\medskip
\small
\[
\begin{array}{c|c}
 \underline{Br(\phi)\ (\textrm{mod}\  \Z^2)}
& \underline{Br(\psi)\ (\textrm{mod}\  \Z^2)}\\[.1in] 

(M-I)^{-1}\left(\!\!\begin{array}{r}-1\\0\end{array}\!\right)=\begin{pmatrix}\frac{3}{4}\\\frac{3}{4}\end{pmatrix} 
& 
\ (M-I)^{-1}\left(\!\!\begin{array}{r}-2\\0\end{array}\!\right)=\begin{pmatrix}\frac{1}{2}\\\frac{1}{2}\end{pmatrix}\\[.3in]
(M-I)^{-1}\left(\!\!\begin{array}{r}2\\0\end{array}\!\right)=\begin{pmatrix}\frac{1}{2}\\\frac{1}{2}\end{pmatrix} & (M-I)^{-1}\left(\!\!\begin{array}{r}3\\0\end{array}\!\right)=\begin{pmatrix}\frac{3}{4}\\\frac{3}{4}\end{pmatrix}
\\[.3in]
(M-I)^{-1}\left(\!\!\begin{array}{r}1\\1\end{array}\!\right)=\begin{pmatrix}\frac{1}{4}\\0\end{pmatrix} &(M-I)^{-1}\left(\!\!\begin{array}{r}1\\2\end{array}\!\right)=\begin{pmatrix}\frac{1}{4}\\\frac{3}{4}
\end{pmatrix}
\end{array}
\]
\vspace{.1in}
\normalsize
There are infinitely many different triangles on the 2-torus having a given set of three points as vertices. 
In case the vertices are the points of  $Br(\phi)$, we write the general edge vectors, up to integer translation in the lift, as 
$(\frac{1}{4}+p,\frac{1}{2}+q) $ and $(\frac{1}{4}+r,\frac{1}{4}+s)$. The area 
has the form
\[
\left|\frac{1}{16}-\frac{1}{8}+\frac{1}{4}(p+s+4ps-q-2r-4rq)\right| \text{.}
\]
Since $p,q,r,s$ are integers, the minimum area is $\frac{1}{16}$ (e.g. $p=q=r=s=0$). 
The argument for $Br(\psi)$ is similar, but the minimum area is $\frac{1}{8}$. Were $\T_\phi$ and $\T_\psi$ homeomorphic, the map $\tau\circ T$ of Corollary~\ref{pisotpart} would preserve the areas of these triangles. We conclude they are not homeomorphic. \qed
\newpage
 \begin{example}\end{example}
 The pair 
 \begin{equation*}
 \phi:=\begin{cases}a\to aabbcac\\ b\to ababccbcc\\c\to aaa^4b^6c^8
 \end{cases}\qquad
  \psi:= 
 \begin{cases}a\to caabbac\\ b\to cababcbcc\\c\to caaa^4b^6c^7
 \end{cases}
 \end{equation*}
 provides an example of substitutions $\phi$ and $\psi$ such that the Perron-Frobenius eigenvalue $\lambda_\phi=\lambda_\psi$ is Pisot but not a unit. The reader can check that the branch locus $Br(\psi)$ 
is colinear, whereas $Br(\phi)$ is not. By Theorem~\ref{main}, $\T_\phi$ and $\T_\psi$ are not homeomorphic.  \qed

\sekshun{Proof of Theorem~\ref{main}} Suppose that $\phi$ is a primitive and
aperiodic substitution with \emph{ language} $$\L_\phi :=
\{w \in \A^*: w \textrm{ is a factor of } \phi^n(i) \textrm{ for some } i \in \A,
n \in \N\}$$ and \emph{ allowed bi-infinite words} $$X_\phi :=
\{ \cdots x_{-1} x_0 x_1 \cdots\, |\  x_n \in \A, x_n  \cdots x_{n+k} \in \L_\phi
\textrm{ for } n \in \Z, k \in \N\}.$$  Suppose 
 that the finite subset $W = \{\varw_1, \ldots, \varw_k\} \subset \L_\phi$ has the
properties that (1) each element of $X_\phi$ can be uniquely factored as a
product of elements of $W$, and (2) for each $\varw_i \in W$, $\phi(\varw_i)$ can be
uniquely factored as a product of elements of
$W$.  We may then define a substitution $\phi': \{1, \ldots, k\} \arrow
\{1, \ldots, k\}^*$ by the rule $\phi'(i) = i_1 \cdots i_\ell$ provided
$\phi(\varw_i) = \varw_{i_1} \cdots \varw_{i_\ell}$.  We will call such a $\phi'$ a \emph{
rewriting of $\phi$}, and the  morphism $\rho: \{1, \ldots, k\}
\arrow \A^*$ given by $\rho(i) = \varw_i$ the associated \emph{ rewriting morphism}. 
Note that $\rho \circ \phi' = \phi \circ \rho$.  It is not hard to see that $\phi'$
is also primitive and aperiodic.

Suppose now that $\phi'$ and $\phi$ are substitutions and $\rho$ is a morphism
with  $\rho \circ \phi' = \phi \circ \rho$. Then the abelianizations
$A, A'$ and $ P$ of $\phi, \phi'$ and $ \rho$, resp., satisfy
$PA' = AP$.  Consequently, if $w$ and $w'$ are positive right 
Perron-Frobenius eigenvectors of $A$ and $A'$, resp., then 
$Pw' = aw$ for some $a > 0$, and $A$ and $A'$ have the same Perron-Frobenius
eigenvalue.  Thus $\phi'$ is Pisot if and only if
$\phi$ is.   In this case, let $V$ and $V'$ be the invariant
Pisot subspaces (with rational bases) corresponding to $A$ and $A'$, resp.  Since $w \in PV' \cap V$, $PV' \cap V$ is a nontrivial rational subspace of $V$.  The characteristic polynomial of $A|_V$ is irreducible, and $PV' \cap V$ is $A$-invariant, which implies that  $PV' \cap V = V$ and $V \subset PV'$.  But $\dim V = d = \dim V'$,    hence $P(V') = V$ and $P|_{V'} : V' \arrow V$ is invertible.     Furthermore, if
$\phi$ and $\phi'$ are Pisot, there is a natural continuous surjection
$\overline{\rho}: \T_{\phi'} \arrow \T_\phi$ that conjugates $ {\Phi}'$ with $ {\Phi}$,
defined as follows:  Given an edge $I = \sigma_1' + x$ in $V'$, with $\rho(i) =
i_1 \cdots i_l$ and $\varv_j := \pi_V e_j$, let
$\overline{\rho}(I)$ denote the finite strand in $V$ defined by
$$\overline{\rho}(I) :=
\{\sigma_{i_1} + Px, 
\sigma_{i_2} + Px + \varv_{i_1}, \ldots, \sigma_{i_\ell} + Px + \varv_{i_1} +
\cdots + \varv_{i_{\ell-1}}\}.$$  If $\gamma' = \{\sigma'_{i_k} + x_k: k \in
\Z\}$ is any strand in $V'$, let $\overline{\rho}(\gamma')$ denote the strand in
$V$ made up of the  finite substrands 
$\overline{\rho}(\sigma_{i_k} + x_k)$ given by  
\[
\overline{\rho}(\gamma') =
\bigcup_{k \in \Z}\overline{\rho}(\sigma'_{i_k} + x_k).
\]  
Observe 
that $\overline{\rho} \circ  {\Phi}' =  {\Phi} \circ
\overline{\rho}$.  Since $P\E^u_{\phi'} = \E^u_\phi$, there is an $R > 0$ such that
$\overline{\rho}(\T_{\phi'}) \subset \F^R_\phi$.  Thus if $\gamma' \in
\T_{\phi'}$, 
\[\overline{\rho}(\gamma') = \overline{\rho}\left(( {\Phi}')^n ( {\Phi}')^{-n} (\gamma')\right) =
 {\Phi}^n(\overline{\rho}\left(( {\Phi}')^{-n} (\gamma')\right) \in  {\Phi}^n \F^R_\phi
\] 
for all $n \in
\N$,  so 
$$\overline{\rho}(\gamma')
\in
\bigcap_{n
\in
\N} {\Phi}^n(\F^R_\phi) = \T_\phi.$$  Continuity of $\overline{\rho}$ is clear, and from
the fact that $\overline{\rho}(\gamma' - t) =
\overline{\rho}(\gamma' ) - at$ and minimality of the tiling flow, it follows that
$\overline{\rho}: \T_{\phi'} \arrow \T_\phi$ is surjective.

Suppose that $V$ and $V'$ are finite dimensional real vector spaces and $\Sigma
\subset V$, $\Sigma' \subset V'$ are lattices. By a map
$T:(V, \Sigma) \arrow (V', \Sigma')$, we will mean a linear transformation $T :
V \arrow V'$ with $T(\Sigma) \subset \Sigma'$.  Two maps $A: (V, \Sigma) \arrow
(V, \Sigma)$ and $A': (V', \Sigma')
\arrow (V', \Sigma')$ are \emph{ shift-equivalent}, $A
\sim_{SE} A'$, provided there are maps $T: (V, \Sigma) \arrow (V',
\Sigma')$ and $S: (V', \Sigma') \arrow (V, \Sigma)$ and natural numbers $m, n$
so that $A^m = ST$, $(A')^n = TS$, $A'T = TA$, and $SA' = AS$;  $S$ and $T$ \emph{
express the shift equivalence}. The relation $\sim_{SE}$ is an equivalence
relation, and if $A
\sim_{SE} A'$, then
$A^k
\sim_{SE}(A')^k$ for all $k \in \N$.  Furthermore, if $\hat{F}_T:
\inv F_A \arrow \inv F_{A'}$ is defined by $$\hat{F}_T (x_1 + \Sigma, x_2 +
\Sigma,
\ldots) = (Tx_1 + \Sigma', Tx_2 + \Sigma', \ldots)$$  (and
$\hat{F}_S, \hat{ F}_{A'} $ and $\hat{ F}_{A}$ are defined similarly), and  if
$\Sigma$ and
$\Sigma'$ are co-compact, then $\hat{F}_T$  and
$\hat{F}_S$ are topological isomorphisms of the solenoids that conjugate the
shifts $\hat{F}_A$ and $\hat{F}_{A'}$.

\begin{lemma} \label{rewriting} Suppose that $\phi$ and $\phi'$ are Pisot
substitutions and that either \\ (a)  $\phi'$ is a rewriting of $\phi$ with
rewriting morphism $\rho$, or \\ (b)  there are morphisms
$\alpha, \beta$ so that $\phi = \alpha \circ \beta$ and $\phi' = \beta
\circ \alpha$. \\ Let
$A, V$ and $ \Sigma$ ($A', V'$ and $ \Sigma'$, resp.) be the
abelianization, invariant Pisot subspace and return lattice for $\phi$
($\phi'$, resp.).  Then, in the case of (a), \\ (1a)
$\overline{\rho}:
\T_{\phi'} \arrow \T_\phi$ is a homeomorphism,  \\ or, in the case of (b),   \\
(1b)
$\overline{\alpha}:
\T_{\phi'} \arrow \T_\phi$ and $\overline{\beta}:
\T_{\phi} \arrow \T_{\phi'} $ are homeomorphisms.  \\Also, \\(2)  the maps
$A|_V: (V, \Sigma) \arrow (V, \Sigma) $ and $A'|_{V'}: (V', \Sigma')
\arrow (V', \Sigma') $  are shift equivalent. \\ Furthermore, if
$\varg_\phi : \T_\phi \arrow \inv F_{A|_V}$ and $\varg_{\phi'} :
\T_{\phi'}
\arrow \inv F_{A'|_{V'}}$  are geometric realizations onto the natural
solenoids, and $T: (V, \Sigma) \arrow (V',
\Sigma')$ and $S: (V', \Sigma') \arrow (V, \Sigma)$ express the shift
equivalence in (2), then there are translations $\tau$ and $\tau'$ on
$\inv F_{A|_V}$ and $\inv F_{A'|_{V'}}$ so that, in the case of (a),  \\ (3a)
$
\varg_\phi
\circ \overline{\rho} = \tau \circ \hat{F}_S \circ \varg_{\phi'}$, and \\
\hspace*{.3in}  $ \varg_{\phi'}
\circ \overline{\rho}^{\,-1} = \tau' \circ \hat{F}_T \circ
\hat{F}^m_{A|_V}\circ
\varg_{\phi}$ for some $m \in \N$, \\ or, in the case of (b),
\\ (3b)  $ \varg_\phi
\circ \overline{\alpha} = \tau \circ \hat{F}_S \circ \varg_{\phi'}$, and \\
\hspace*{.3in}  $ \varg_{\phi'}
\circ \overline{\beta} = \tau' \circ \hat{F}_T \circ
\varg_{\phi}$.
\end{lemma}

\begin{proof} Suppose that $\phi'$ is a rewriting of $\phi$ with rewriting morphism
$\rho$ having abelianization $P$, and suppose that $\gamma'$,
$\gamma'' \in \T_{\phi'}$ are such that $\overline{\rho} (\gamma') =
\overline{\rho} (\gamma'')$.  Let $w'$ and $w''$ be the words of $X_{\phi'}$ spelled
out by  $\gamma'$ and 
$\gamma'' $, where some location is chosen for the decimal point indicating the
location of the $0^{th}$ letter.  Then $\rho(\varw')$ and 
 $\rho(\varw'')$ agree, up to shift, so that $w'$ and $w''$ must agree, up to shift,
by the uniqueness of factorization in $X_\phi$ in the definition of rewriting. 
It follows that $\gamma' = \gamma''-t$ for some $t $.  But then 
$\overline{\rho} (\gamma'') = \overline{\rho} (\gamma') = \overline{\rho} (\gamma'' - t)
=\overline{\rho} (\gamma'') - at$, hence $t = 0$ and $\gamma' = \gamma''$.  That is,
$\overline{\rho}$ is a homeomorphism.

In case $\phi = \alpha \circ \beta$ and $\phi' = \beta \circ \alpha$, we have that
$\overline{\alpha} \circ \overline{\beta} =  {\Phi}$ and   $
 \overline{\beta} \circ \overline{\alpha}=  {\Phi}'$, so both $\overline{\alpha}$ and
$\overline{\beta}$ must be homeomorphisms.

Again, suppose that $\phi'$ is a rewriting of $\phi$.  If $v  \in
\Theta'(i)$ is a return vector for $\phi'$, say $\gamma' \in
\T_{\phi'}$ has edges $I$ and $I + v$ of type $i$, then $\overline{\rho}
\gamma' \in \T_\phi$ has edges $J$ and $J + Pv$ of type $i_1$, where
$\rho(i) = i_1 \cdots$.  Thus $Pv \in \Theta(i_1)$ is a return vector for
$\phi$.  It follows that $ P \Sigma' \subset \Sigma$.

While it is not necessarily the case that $P\Sigma' = \Sigma$, we shall see that
$A^m \Sigma \subset P \Sigma'$ for some $m \in \N$.  To this  end, suppose that
$I$ is an edge in $V$ of type $i$, and $I$, $I+v$ are edges of $\gamma \in
\T_\phi$ with $v \neq 0$.  Let $\gamma' :=
\overline{\rho}^{\,-1}(\gamma)$, and for each $n \in \N$, let $\delta^n$ and $\eta^n$ be
maximal substrands of $( {\Phi}')^n(\gamma')$  with the property that $\overline{\rho}(\delta^n) \subset
 {\Phi}^n(I)$ and  $\overline{\rho}(\eta^n) \subset
 {\Phi}^n(I + v)$.  If $x^n, y^n \in \L_{\phi'}$ are the words corresponding to
$\delta^n$, $\eta^n$, resp., then we have the factorizations $p^n
\rho(x^n)s^n$ and  $q^n \rho(y^n)r^n$ of $\phi^n(i)$ in which the words $p^n$,
$s^n$, $r^n$, $q^n \in \L_\phi$ are of bounded length. Putting  a decimal point
(to mark the position of the $0^{th}$ coordinate) on the left of $\phi^n(i)$, we
may choose $n_k \arrow \infty$ and $m_k \in \N$, with $m_k$ on the order of
$\frac{n_k}{2}$, so that $\sigma^{m_k}(\phi^{n_k}(i))$ converges to a bi-infinite
word $w \in X_\phi$ (here $\sigma$ is the shift that takes $.abc \cdots$ to 
$a.bc \cdots$).  Let $i^k$ be the letter (in $\A = \A_\phi$) immediately to the right of the decimal point in $\sigma^{m_k}(\phi^{n_k}(i))$, and let $x_0^{n_k}$ and $y_0^{n_k}$ denote the letters (in $\A' = \A_{\phi'}$)  of  $x^{n_k}$ and $y^{n_k}$, resp.,  so that $\rho( x_0^{n_k}) $ and $\rho( y_0^{n_k})$ contain (the identified  occurrence of) $i^k$.   Putting a decimal point immediately to the left of  $x_0^{n_k}$ and $y_0^{n_k}$  in  $x^{n_k}$ and $y^{n_k}$, resp., we may choose a subsequence $n_{k_\ell}$
so that  $x^{n_{k_\ell}}$ and $y^{n_{k_\ell}}$ converge to bi-infinite words, say $x$ and $y$, resp.,  in $X_{\phi'}$.  Now $\rho (x)$ and $\rho (y)$ give factorings of $\varw \in X_\phi$ by the words $\varw_i = \rho(j)$, for $j \in \A'$.  By uniqueness,       $x = y$.  Thus for large $\ell$, $x_0^{n_{k_\ell}} = y_0^{n_{k_\ell}}$.  Pick such an $\ell$, let $m = {n_{k_\ell}}$, $j = x_0^m = y_0^m$, and $J,  J + \varv'$ be the edges of the substrands $\delta^m$, $\eta^m$ of $({\Phi}')^m(\gamma')$ corresponding to  $x_0^m $ and $ y_0^m$, resp. Then $P (\varv' ) = A^m \varv$.  That is, given $v \in \Theta(i)$, a return vector for $\phi$, there is $\varv' \in \Theta'(j)$, a return vector for $\phi'$,   with $\varv = (A|_V)^{-m}P(\varv') = P(A'|_{V'})^{-m}(\varv')$.

Now suppose that $\varv_1, \ldots, \varv_d$ is a basis for $\Sigma$, the return lattice for $\phi$.  Each $\varv_j$ is then an integer combination of elements of $\bigcup_{k \in \Z} (A|_V)^k \Theta(i)$.  The preceding argument shows that there is an $m \in \N$ so that $A^m \varv_j$ is in $P \Sigma'$ for all $j = 1, \ldots, d$.  That is, $A^m \Sigma \subset P\Sigma'$. Let $S:= P|_{V'}$ and $T:= (P|_{V'})^{-1}(A|_V)^m$
(recall that $P|_{V'}$ is invertible because $\phi$ and $\phi'$ are Pisot).
We have $S: (V', \Sigma') \arrow (V, \Sigma)$, $T: (V, \Sigma) \arrow (V',
\Sigma')$, $A|_V S = S A'|_{V'}$, $T A|_V = A'|_{V'}T $, $ST = (A|_V)^m$, and $TS = (A'|_{V'})^m$.

Still in the rewriting case, suppose that $\varg_\phi: \T_\phi \arrow \inv F_{A|_V}$ and $\varg_{\phi'}: \T_{\phi'} \arrow \inv F_{A'|_{V'}}$ are geometric realizations onto the natural solenoids defined by the selection of $\{u_i\}$ in $\Gamma/\Sigma$ and $\{u'_i\}$ in $\Gamma'/\Sigma'$, resp. (see the definition of geometric realization in \S\ref{geom}).  Let $\{\varw_{ij}\} \subset \Gamma/\Sigma$ and $\{w'_{ij}\} \subset \Gamma'/\Sigma'$ denote the corresponding transition vectors (so $\varw_{ij} = u_j - u_i$  and $w'_{ij} = u'_j - u'_i$, for all $i, j$).  Let $\gamma' \in \T_{\phi'}$ have edge $I$ of type $i$.  Then $\varg_{\phi'}(\gamma')$ has $0^{th}$coordinate $$(\varg_{\phi'}(\gamma'))_0 = (\min I + \Sigma') - u'_i $$ in  $V'/\Sigma'.$  The strand $\overline{\rho}(\gamma')$ in $\T_\phi$ has an edge $J$ of type $i_1$, where $\rho(i) = i_1 \cdots$, with $\min J = P \min I$.  Thus $$(\varg_\phi(\overline{\rho}(\gamma')))_0 = ((P \min I )+ \Sigma) -u_{i_1}.$$
We have 
\begin{eqnarray*} F_P((\varg_{\phi'}(\gamma'))_0)& =  & \left((P \min I\right) + \Sigma) - F_P u'_i \\ &= & \left(\varg_\phi(\overline{\rho}(\gamma'))\right)_0 + u_{i_1} - F_P u'_i. 
\end{eqnarray*}  We claim that $u_{i_1} - F_P u'_i $ is independent of $i \in \A'$.  Indeed, 
 \begin{eqnarray*} (u_{i_1} - F_P u'_i ) -  (u_{j_1} - F_P u'_j) &=& (u_{i_1} - u_{j_1})  - (F_P (u'_i -   u'_j) )\\ &=& \varw_{j_1 i_1} - F_P(w'_{ji}).\end{eqnarray*}   But if $\gamma' \in \T_{\phi'}$ has edge $I$ of type $i$ and edge $J$ of type $j$, so that $(\min J - \min I) + \Sigma = w'_{ji}$, then $\overline{\rho}(\gamma')$ has corresponding edges of type $i_1$ and $j_1$ with initial vertices differing by $P( \min J) - P(\min I)$.  That is, $\varw_{j_1 i_1} = F_P(\varw'_{ji})$.  It follows that if $\tau_1$ is translation by $F_P u'_i - u_{i_1}$, then $$(\varg_\phi \circ \rho)_0 = \tau_1 \circ  F_{P|_{V'}} \circ (\varg_{\phi'})_0.$$ Similarly,  $k^{th}$ coordinates satisfy $$(\varg_\phi \circ \rho)_k = \tau_1 \circ  F_{P|_{V'}} \circ (\varg_{\phi'})_k.$$
 
Note that \begin{eqnarray*}F_{A|_V}(u_{i_1} - F_{P|_{V'}}  u'_i )& =& F_{P|_{V'}} u_{i_1} - F_{P|_{V'}} F_{A|_{V'}} u'_i \\& =& u_{i_1} - F_{P|_{V'}}  u'_i \end{eqnarray*} from the normalization requirement on the $\{u_i\}, \{u'_i\}$,  so translation by this same element in each coordinate in $\inv F_{A|_V}$  defines a translation $\tau$ on $\inv F_{A|_V}$.
We have $$\varg_\phi \circ \overline{\rho} = \tau \circ \hat{F}_{P|_{V'}} \circ \varg_{\phi'}.$$  Thus with $S = P|_{V'}$ and $T = (P|_{V'})^{-1} \circ (A|_{V})^m$, we have conclusion (3a).

The proofs of (2), in case $\phi = \alpha \circ \beta$ and $\phi' = \beta \circ \alpha$, and (3b) are similar (although more straightforward).  \end{proof}

\

With the notation of Lemma \ref{rewriting}, and with $C_\phi$ and
$C_{\phi'}$ the collections of special tilings for $\phi$ and $\phi'$
resp., note that $\overline{\rho}$, $\overline{\alpha}$, and $\overline{\beta}$ determine
bijections between $C_\phi$ and
$C_{\phi'}$, and $\hat{F}_{A|_V}$ maps $\varg_\phi(C_\phi)$ bijectively onto
itself.  This yields the following consequence of Lemma \ref{rewriting}.

\begin{corollary} \label{corollary} If $\phi$ and $\phi'$ are Pisot substitutions
with one a rewriting of the other, or for which there are morphisms
$\alpha, \beta$ with $\phi= \alpha \circ \beta$ and $\phi' = \beta \circ
\alpha$, then $A|_V: (V, \Sigma) \arrow (V, \Sigma) $ and $A'|_{V'}: (V',
\Sigma')
\arrow (V', \Sigma') $  are shift equivalent, expressed by maps  $S: (V',
\Sigma') \arrow (V, \Sigma)$ and $T: (V, \Sigma) \arrow (V',
\Sigma')$, with $$\tau \circ F_T\circ \pi_0 \circ \varg_\phi(C_\phi) =
\pi_0 \circ \varg_{\phi'}(C_{\phi'})$$ and $$\tau' \circ F_S \circ \pi_0 \circ
\varg_{\phi'}(C_{\phi'}) =
\pi_0 \circ \varg_{\phi}(C_{\phi}).$$
\end{corollary}
Let us say that the substitutions $\phi$ and $\psi$ are in the same 
rewriting class if  there are substitutions  $\phi_0=\phi,\ \phi_1,\ \dots, \phi_n=\psi$ with the 
property that for each $i=0,\dots,n-1$,  one of $\phi_i$ and $\phi_{i+1}$
is a rewriting of the other. The proof of the rigidity result in \cite{bargeswanson2} (see Theorem~\ref{rigidity} in the next section) establishes the following:

\begin{lemma}\label{rigiditythm}  Suppose that  $\phi$ and
$\psi$ are substitutions such that $\T_\phi$ is homeomorphic with
$\T_\psi$.  Then there are $m, n \in \N$ and substitutions  $\phi'$ and
$\psi'$ such that (1) $\phi'$ and $\phi^m$ are in the same rewriting class, as are
$\psi'$ and $\psi^n$, and (2) there are morphisms
$\alpha$ and $\beta$ such that $\phi'  = \alpha \circ \beta$ and $\psi' = \beta
\circ \alpha$.
\end{lemma}

Proof:  See \cite[Theorem 2.1]{bargeswanson2}.  \qed

\medskip

To prove Theorem \ref{main}, let $\phi$ and $\psi$ be Pisot of degree
$d$, with geometric realizations $\varg_\phi$ and $\varg_\psi$ of
$\T_\phi$ and $\T_\psi$.  Fix isomorphisms $L_\phi : (V_{\phi},
\Sigma_\phi) \arrow (\R^d,
\Z^d)$ and  $L_\psi : (V_{\psi}, \Sigma_\psi) \arrow (\R^d, \Z^d)$, and let 
$Br(\phi) = F_{L_\phi}\circ \pi_0 (\varg_\phi(C_\phi))$ and $Br(\psi) =
F_{L_\psi}\circ \pi_0 (\varg_\psi(C_\psi))$ be the corresponding branch loci.  Let $m,
n, \phi'$ and $\psi'$ be as in Lemma~\ref{rigiditythm}.  As $\T_\phi =
\T_{\phi^m}$ and   $\T_\psi =
\T_{\psi^n}$, we may take $\varg_{\phi^m} = \varg_\phi$, $\varg_{\psi^n} = \varg_\psi$,
$L_{\phi^m} = L_\phi$ and $L_{\psi^n} = L_\psi$, so that 
$Br(\phi^m)= Br(\phi)$  and $Br(\psi^m) = Br(\psi) $.  Apply Corollary
\ref{corollary}  repeatedly to get translations $\eta$ and
$\hat{\eta}$ (this last determined by $\eta$ and $L_\psi$), and  a map $T_1 : (V_\phi,
\Sigma_\phi) \arrow (V_\psi, \Sigma_\psi)$ for which there are $k, \ell \in \N$ and a map 
$S_1 :(V_\psi,\Sigma_\psi) \arrow (V_\phi, \Sigma_\phi)$ such that $$S_1 T_1 = (A_\phi^m|_{V_\phi})^k, \quad T_1S_1
= (A_\psi^n|_{V_\psi})^\ell$$ and 
\begin{eqnarray*} Br(\psi)& = & F_{L_\psi}\left((\varg_{\psi^n}(C_{\psi^n}))_0\right) \\ & = & F_{L_\psi}(\eta \circ
F_{T_1}\left(\varg_{\phi^m}(C_{\phi^m}))_0\right)\\ & = & 
F_{L_\psi}\left(\eta \circ
F_{T_1}\left(L^{-1}_\phi(Br(\phi)\right)\right)\\ & = & \hat{\eta} \circ F_{L_\psi T_1 L^{-1}_\phi}
.\end{eqnarray*}

If $M_\phi = L_\phi A_\phi L^{-1}_\phi$, $M_\psi = L_\psi A_\psi L^{-1}_\psi$, $T =
L_\psi T_1 L^{-1}_\phi$, and $S = L_\phi S_1 L^{-1}_\psi$ (expressed as matrices in the standard basis on $\R^d$), $m_0 = m^k$ and $\tau_0 = \hat{\eta}$,
then $Br(\psi) =
\tau_0 \circ F_T(Br(\phi))$ with $ST = (M_\phi)^{m_0}$. Similarly, if $m_1=n^\ell$, there is a translation $\tau_1$ with $Br(\phi)=\tau_1\circ F_S\left((Br(\psi)\right)$ and $TS=(M_\psi)^{m_1}$.  \qed


\sekshun{The Pisot Part of the Augmented Dimension Group}\label{pisotpartsection}
In this section we give a cohomological interpretation of Theorem~\ref{main}. For a Pisot substitution $\phi$ of degree $d$, let $L$ and $M$ be as in the definition of $Br(\phi)$. The induced 
$F_M:H^1(\cT^d;\mathbb R)\hookleftarrow$ is then a vector space isomorphism with simple eigenvalue $\lambda=\lambda_\phi$ and a $(d-1)$-dimensional invariant subspace complementary to the eigenspace of $\lambda$. This $(d-1)$-dimensional space splits 
$H^1(\cT^d;\R)$ into two closed invariant half-spaces, one of which contains the cocycle dual to $(F_L)_*([\omega])$, where $[\omega]$ is the 1-cycle generated by the positive eigenvector $\omega$ of $A=A_\phi$. This half-space is called the \emph{positive cone} in $H^1(\cT^d;\R)$ and is denoted 
 by $H^1(\cT^d;\R)^+$. The inclusion $\imath:\Z\to\R$ induces $\imath^*:H^1(\cT^d)\to H^1(\cT^d;\R)$ (an unspecified coefficient ring is understood to be $\Z$), and we define the positive cone in $H^1(\cT^d)$ to be $H^1(\cT^d)^+:= (\imath^*)^{-1}(H^1(\cT;\R)^+)$. Likewise, the positive cone in $H^1(\cT^d,Br(\phi))$ is
\[
H^1(\cT^d,Br(\phi))^+:=(\jmath^*)^{-1}\left(H^1(\cT^d)^+\right),
\]
$\jmath^*:H^1(\cT^d,Br(\phi))\to H^1(\cT^d)$ the natural homomorphism.

We have \[ F_M^*:\left(H^1(\cT^d,Br(\phi)),(H^1(\cT^d,Br(\phi))^+\right)\to\left(H^1(\cT^d,Br(\phi)),(H^1(\cT^d,Br(\phi))^+\right),\] and we define $PDG(\phi):=\inv \left(F_B^*:\, H^1\left(\cT^d,Br(\phi)\right) \to H^1\left(\cT^d,Br(\phi)\right)\right)$, with positive cone 
\[
PDG(\phi)^+:=\left\{[(k,g)]\in PDG(\phi)\,:\ g\in H^1\left(\cT^d,Br(\phi)\right)^+\right\} 
\]
The \emph{shift isomorphism} on $\left(PDG(\phi),PDG(\phi)^+\right)$ is the ordered isomorphism given by $[(k,g)]\mapsto \left[\left(k,F_M^*(g)\right)\right]$.

A homeomorphism of tiling spaces, $f:\T_\phi\to\T_\psi$, is \emph{orientation preserving (reversing)} provided it takes the positive flow direction in $\T_\phi$ to the positive (negative) flow direction  in $\T_\psi$; that is, the function $s:\,\T_\phi\times\R^+\to \R$ defined by $f(\gamma - t)=f(\gamma)-s(\gamma,t)$ is positive (negative).
\begin{theorem}\label{pisotpart}
If $\phi$ and $\psi$ are Pisot substitutions, and the tiling spaces $\T_\phi$ and $\T_\psi$ are orientation preserving (reversing) homeomorphic, then there is an ordered isomorphism between $\ds \left(PDG(\phi),PDG(\phi)^+\right)$ and $\ds \left(PDG(\psi),PDG(\psi)^+\right)\  \left({\mathrm resp.}\ \ds \left(PDG(\psi),-PDG(\psi)^+\right)\right)$ that conjugates some positive powers of the shift isomorphisms. 
\end{theorem}
\begin{proof}
Let $T,\ S,\ \tau_0,\ \tau_1,\ m_0,\ {\rm and}\ m_1$ be as in Theorem~\ref{main}. The commuting diagram

\setlength{\unitlength}{1in}

\begin{picture}(6,2)
\put(.5,1.5){$\ds H^1\left(\cT^d,Br(\phi)\right) $}
\put(1.7,1.55){\vector(1,0){.5}}
\put(1.85,1.7){$\ds F_{M_\phi^{m_0}}^*$}
\put(2.4,1.5){$\ds H^1\left(\cT^d,Br(\phi)\right)$}
\put(2.9,1.3){\vector(1,-1){.3}}\put(3.2,1.15){$\ds (\tau_1\circ F_S)^*$}
\put(1,1.3){\vector(1,-1){.3}}\put(.35,1.1){$\ds (\tau_1\circ F_S)^*$}
\put(1.8,1){\vector(2,1){.6}}\put(2.2,1.05){$\ds (\tau_0\circ F_T)^*$}
\put(3.0,.7){$\ds H^1\left(\cT^d,Br(\psi)\right)$}
\put(1.2,.7){$\ds H^1\left(\cT^d,Br(\psi)\right)$}
\put(2.4,.75){\vector(1,0){.5}} 
\put(2.55,.55){$\ds F_{M_\psi^{m_1}}^*$}
\end{picture}
\newline
induces an isomorphism $\ds \left[(k,g))\right]\to \left[(k,(\tau_1\circ F_S)^*(g)\right]$ that conjugates the ${m_1}^{\rm th}$ power of the shift on $PDG(\psi)$ with the ${m_0}^{\rm th}$ power of the shift on $PDG(\phi)$. Since the order is dynamically defined, this conjugacy alone guarantees that order is either preserved or reversed. If $\T_\phi$ and $\T_\psi$ are orientation preserving homeomorphic, then  $S$ and $T$ take $L_\psi(\omega_\psi)$ and $L_\phi(\omega_\phi)$ to positive multiples of $L_\phi(\omega_\phi)$ and $L_\psi(\omega_\psi)$, resp. (here, $L_\phi$ and $L_\psi$ denote the $L$ in the definition of $Br(\phi)$  and $Br(\psi)$, and $\omega_\phi$, $\omega_\psi$ are positive eigenvectors of $A_\phi$, $A_\psi$, resp.). Thus 
$(\ds \tau_1\circ F_S)^*\left(H^1(\cT^d,Br(\phi))^+\right)=H^1(\cT^d,Br(\psi))^+$.

If $\T_\phi$ and $\T_\psi$ are orientation reversing homeomorphic, then $\T_\phi$ and $\T_{\overline\psi}$ are orientation preserving homeomorphic, $\overline\psi$ being the reverse of $\psi$ (if $\psi(i)=i_1\cdots i_k$ then 
$\overline\psi(i)=i_k\cdots i_1$). The above applied to $(PDG(\phi),PDG(\phi)^+)$ and  $(PDG(\psi),PDG(\psi)^+)$ together with the isomorphism $(PDG(\overline \psi ),PDG(\overline\psi)^+)\to  (PDG(\psi),-PDG(\psi)^+)$ induced by $F_{-{\rm Id}}^*$ yields the conclusion of the theorem in the orientation reversing case.
 \end{proof}
\begin{example}\end{example}
Consider $\ds \phi:\begin{cases}a\quad \arrow & ababba\\ b\quad \arrow & aabbba \end{cases}$ and  $\psi:\begin{cases}a\quad \arrow & a^{47}ab^{18}bb^{29} \\ b\quad \arrow & a^{47}bb^{18}ab^{29} \end{cases}$.

One computes: 
\[
\begin{array}{ll}\left(PDG(\phi),PDG(\phi)^+\right)\cong & \left(DG\begin{pmatrix}6&0\\1&1\end{pmatrix},DG\begin{pmatrix}6&0\\1&1\end{pmatrix}^+\right),\quad\text{and} \\[.3in]
\left(PDG(\psi),PDG(\psi)^+\right)\cong & \left(DG\begin{pmatrix}96&0\\19&1\end{pmatrix},DG\begin{pmatrix}96&0\\19&1\end{pmatrix}^+\right).\end{array}
\]
These groups are order isomorphic but, because no power of $96$ is a power of $6$, the shift isomorphisms do not have conjugate powers. By Theorem~\ref{pisotpart}, $\T_\phi$ and $\T_\psi$ are not homeomorphic.   \qed

We now relate the pair $(PDG(\phi),PDG(\phi)^+)$ to the \emph{ordered augmented cohomology}  of the tiling space 
$\T_\phi$. We begin with a description of the \emph{augmented tiling space} (for more detail, and a precise description as an inverse limit, see \cite{bargesmith}). 

Let $\phi$  be any primitive aperiodic substitution with tiling space $\T_\phi$. Let 
$\{\gamma_1^f,\dots, \gamma_{n_f}^f\}$ be a collection of forward asymptotic special tilings, exactly one chosen from each forward asymptotic equivalence class, and $\{\gamma_1^b,\dots,\gamma_{n_b}^b\}$ a collection of backward asymptotic special tilings, one from each backward asymptotic equivalence class.  Let $R_j^f$, $j=1,\dots,n_f$,  and $R_i^b$, $i=1,\dots n_b$,  be the rays 

\[
R_j^f:= \left\{\{j\}\times \{\gamma_j^f-t\}\,:\ t\ge 0\right\}\quad\text{and}\quad R_i^b:= \left\{\{i\}\times \{\gamma_i^b-t\}\,:\ t\ge 0\right\}. 
\] 

The augmented tiling space $\widetilde \T_\phi$ is defined to be the union 
\[
\widetilde \T_\phi:= \T_\phi\ \   \cup\  \ \left.\left(\bigcup_{j=1}^{n_f}R_j^f\,\cup\,\bigcup_{i=1}^{n_b} R_i^b\right) \ \right/ \, 
 \left\{
\{j\}\times \{\gamma_j^f\},\ \{i\}\times \{\gamma_i^b \}
\right\},
\]
in which all of the endpoints of the rays have been identified to a single branch point.  The metric on $\T_\phi$ is extended to $\widetilde\T_\phi$ in such a way that 
\[
d\left(\gamma_j^f-t,\{j\}\times \{\gamma_j^f-t\}\right)=\frac{1}{1+t}=d\left(\gamma_i^b+t,\{i\}\times \{\gamma_i^b+t\}\right),
\]
for $t\ge 0$, making the ray $R_j^f$ asymptotic to the forward orbit of $\gamma_j^f$ and the ray $R_i^b$ asymptotic to the backward orbit of $\gamma_i^b$.

The homeomorphism $\Phi$ extends to $\widetilde \Phi:\widetilde \T_\phi\hookleftarrow$ with $\ds \widetilde \Phi\left(\{j\}\times \{\gamma_j^f-t\}\right)=\{j'\}\times \{\gamma_{j'}^f-\lambda t\}$, where $j'$ is such that  $\Phi(\gamma_j^f)$ is forward asymptotic to $\gamma_{j'}^f$, and similarly for the rays $R_i^b$. We order $H^1(\widetilde \T_\phi)$  in such a way that $\ds \widetilde \Phi^*:\left(H^1(\widetilde \T_\phi),\, H^1(\widetilde \T_\phi)^+\right)\hookleftarrow$ is an order isomorphism. First note that  $\ds \Phi^*\,:\, H^1(\T_\phi;\R)\hookleftarrow$ is a vector space isomorphism with simple eigenvalue $\lambda=\lambda_\phi$ and codimension one invariant subspace complementary to 
the eigenspace of $\lambda$.  This codimension one space splits $H^1(\T_\phi)$ into two closed 
half-spaces, one of which corresponds
to the positive direction of the flow\footnote{Specifically, there is an orientation preserving map of $\T_\phi$ onto the circle $\mathbb T^1$. Let $H^1(\T_\phi;\R)^+$ contain the positive half-line in $H^1(\T^1;\R)$ pulled back to $H^1(\T_\phi;\R)$. }; call that space $H^1(\T_\phi;\R)^+$. Let 
$H^1(\widetilde \T_\phi)^+=(\imath_*\jmath^*)^{-1}\left(H^1(\T_\phi;\R)^+\right)$, where 
$\imath_*:H^1(\T_\phi)\to H^1(\T_\phi;\R)$ is induced by $\imath:\Z\to\R$ and 
$\jmath^*:H^1(\widetilde \T_\phi)\to H^1(\T_\phi)$ is induced by the inclusion $j:\T_\phi\to \widetilde \T_\phi$. The pair $\ds \left(H^1(\widetilde \T_\phi),\,H^1(\widetilde \T_\phi)^+\right)$ is the \emph{ordered augmented cohomology} of $\T_\phi$. Call the isomorphism $\ds \widetilde \Phi^*:\left(H^1(\widetilde \T_\phi),\, H^1(\widetilde T_\phi)^+\right)\hookleftarrow$ the {\bf shift} isomorphism.

The rigidity result of \cite{bargeswanson2} is 
\begin{theorem}\label{rigidity}
If $f:\T_\phi\to \T_\psi$ is a homeomorphism of substitution tiling spaces, then there are $m,n\in{\mathbb N}$ and  $h:\T_\phi\to\T_\psi$ isotopic to $f$ so that $h\circ\Phi^m=\Psi^n\circ h$.
\end{theorem}
\begin{corollary}\label{augmented}
 If $\T_\phi$ and $\T_\psi$ are orientation preserving (reversing) homeomorphic, then the ordered cohomologies $\ds \left(H^1(\widetilde\T_\phi),H^1(\widetilde\T_\phi)^+\right)$ and  $\ds \left(H^1(\widetilde\T_\psi),H^1(\widetilde\T_\psi)^+\right)$
\big(resp., $\ds \left(H^1(\widetilde\T_\psi),-H^1(\widetilde\T_\psi)^+ \right)$\big) are isomorphic by an isomorphism that conjugates some positive powers of the shift isomorphisms. 
\end{corollary}
\begin{proof}
The homeomorphism $h$ of Theorem~\ref{rigidity} extends to a homeomorphism 
$\widetilde h: \widetilde \T_\phi\to\widetilde \T_\psi$ that conjugates $\widetilde \Phi^m$ with $\widetilde \Psi^n$.   
\end{proof}
The corollary ---  without the conjugacy in the conclusion --- appears in \cite{bargesmith}. Also, in that paper, a nonnegative integer matrix $\widetilde A^t$ is constructed with ordered dimension group $\ds \left(DG(\widetilde A^t),DG(\widetilde A^t)^+\right)$ isomorphic to $\ds \left(H^1(\widetilde \T_\phi),H^1(\widetilde \T_\phi)^+\right)$.

Viewing $\ds \left(PDG(\phi),PDG(\phi)^+\right)$ as the Pisot part of the augmented cohomology is justified by our  final theorem. A formula for  the Pisot part as the  dimension group of an integer matrix will arise in the proof.
 
\begin{theorem}\label{pisotcohomology}
If $\phi$ is a Pisot substitution, then there is an ordered embedding of $\ds \left(PDG(\phi),PDG(\phi)^+\right)$ into  $\ds \left(H^1(\widetilde \T_\phi),H^1(\widetilde \T_\phi)^+\right)$ that commutes with  a positive power of the shift isomorphisms.
\end{theorem}

\begin{example}\end{example}
   If $\phi$ is the Thue-Morse substitution ($1\to 12,\ 2\to 21$), then 
\[ \left(PDG(\phi),PDG(\phi)^+\right) = \left(\Z\left[\small\frac{1}{2}\normalsize\right],\, \left\{x\in\Z\left[\small\frac{1}{2}\normalsize\right],\ x\ge 0\right\}\right),\] 
with shift isomorphism $x\mapsto 2x$, and  
\[
\left(H^1(\widetilde \T_\phi),H^1(\widetilde \T_\phi)^+\right)=\left(\Z\left[\small\frac{1}{2}\normalsize\right]\oplus \Z^4,\ \left\{(x,y)\in \Z\text{\tiny$\left[\frac{1}{2}\right]$\normalsize}\oplus \Z^4:\ x\ge 0\right\}\right),
\]
with shift isomorphism $(x,y)\mapsto (2x,y)$. 
\qed
\begin{proof}{(Theorem~\ref{pisotcohomology})} 
Let $\phi$ be a Pisot substitution. To simplify, we pass to a \emph{prepared} substitution (see \cite{bargesmith}).
We may suppose all special asymptotic tilings in $\T_\phi$ are fixed by $\phi^n$, for some $n\in\mathbb N$, and   $\phi^n$ has a fixed bi-infinite word $\cdots b.a\cdots$. So we have $\phi^n(a)=a\cdots$, $\phi^n(b)=\cdots b$
and $ba$ occurs in $\phi^k(a)$, for some $k\in\mathbb N$. 

We rewrite $\phi^n$ with stopping rule $b$ and starting rule $a$ as follows. Let $\mathcal U=\{u_1,\dots, u_q\}$ be the finite collection of words with the properties
 
\begin{itemize}
\item[({\emph i})] each $u_i$ occurs as a factor of the infinite word $a\phi^n(a)\phi^{2n}(a)\cdots $; 
\item[({\emph ii})] each  $u_i$ has the form $a\cdots b$.
\item[({\emph iii})] $ba$ is not a factor of any $u_i,\ i=1,\dots,q$. 
\end{itemize}   
Then each word $\phi^n(u_i)$ factors uniquely in the form  $u_{i_1}\cdots u_{i_p}$. From this, define $\psi$ by $\psi(i)=i_1\cdots i_{p}$, $i=1,\dots,q$. The substitution $\psi$ is primitive and  aperiodic Pisot ($\lambda_\psi=\lambda_\phi^n$) and \emph{strictly proper}: there are $r,s\in\{1,2,\dots, q\}$ so that 
$\psi(i)=r\cdots s$ for $i=1,\dots, q$. Also, the special asymptotic tilings in $\T_\psi$ are all fixed by $\Psi$. Then $\T_{\phi^n}=\T_\phi$ and $\T_\psi$ are orientation preserving homeomorphic (this is a special case of 
Lemma~\ref{rewriting}) by a homeomorphism that conjugates $\Phi^n$ with $\Psi$. Hence,     
$\ds \left(H^1(\widetilde \T_\phi),H^1(\widetilde \T_\phi)^+\right)$ and $\ds \left(H^1(\widetilde \T_\psi),H^1(\widetilde \T_\psi)^+\right)$ are isomorphic by an isomorphism that conjugates $(\Phi^n)^*$ with $\Psi^*$ (Corollary~\ref{augmented}) and $\ds \left(PDG(\phi),PDG(\phi)^+\right)= \left(PDG(\phi^n),PDG(\phi^n)^+\right)$ is isomorphic to $\ds \left(PDG(\psi),PDG(\psi)^+\right)$ (Theorem~\ref{pisotpart}). The isomorphism conjugates the $(\ell n)^{\rm th}$ power of the shift on  $\left(PDG(\phi),PDG(\phi)^+\right)$ to the $\ell^{\rm th}$ power of the shift on $\left(PDG(\psi),PDG(\psi)^+\right)$ for some $\ell\in\mathbb N$. Thus once we prove the theorem for such substitutions $\psi$, we have proved it for all Pisot substitutions.   

Let $n_f$ denote the number of equivalence classes of forward asymptotic special tilings in $\T_\psi$, and $n_b$ the number of special backward classes. As in  \cite{bargesmith}, we may select special asymptotic tilings $\gamma_j^f$, $j=1,\dots, n_f$, one from each forward class, and $\gamma_i^b$, $i=1,\dots, n_b$, one from each backward class, with the properties 
\begin{itemize}
\item [(i)] if $J$ is an edge of $\gamma_j^f$ (resp., $I$ an edge of $\gamma_i^b$) that 
meets $\E^s$, then $J$  (resp., $I$) meets $\E^s$ in its interior;
\item [(ii)] if $k$ is the type of $J$ (resp., of $I$), then $\psi(k)=p_j^fks_j^f$ (resp., $\psi(k)=p_i^bks_i^b$) with $p_j^f$ and $s_j^f$ (resp., $p_i^b$ and $s_i^b$) nonempty and 
\[ \varw_j^f:=\cdots \psi^2(p_j^f)\psi(p_j^f)p_j^f\psi(s_j^f)\psi^2(s_j^f)\cdots\ \text{(resp.,}\  \varw_i^b:=\cdots\psi(p_i^b)p_i^b\psi(s_i^b)\cdots )\] is
the bi-infinite word corresponding to $\gamma_j^f$ (resp., $\gamma_i^b$). 
\end{itemize}
In particular, and this will be important later,  
\[A(\min I)+\mathrm{pr}_V\llbracket p^b_i\rrbracket=\min I \quad \text{and}\quad A(\max J)-\mathrm{pr}_V\llbracket p_j^f\rrbracket=\max J, \]
where $\llbracket u\rrbracket$ denotes the abelianization of the word $u$.


Let $E$ denote the $n\times (n_f+n_b-1)$ matrix with $ij^{\rm th}$ entry 
\[
E_{ij}= 
\begin{cases}
\text{number of occurrences of $i$ in $p_1^b s_j^f,\ \text{if}\ 1\le j\le n_f$} 
\\[.1in]
\text{number of occurrences of $i$ in $p^b_{j-n_f+1}s^f_1,\ \text{if}\ n_f<j\le n_f+n_b-1$.}
\end{cases}
\]
Let $A=A_\psi$ be the abelianization of $\psi$, and let $I$ be the $(n_f+n_b-1)\times (n_f+n_b-1)$ identity matrix.  
The \emph{augmented  matrix} for $\psi$ is 
 \[
\widetilde A=\begin{pmatrix}  A & E\\O & I\end{pmatrix}, 
\]
of size $\n=n+n_f+n_b-1$. The \emph{augmented dimension group} for $\psi$ is the pair
\[
\left( DG(\widetilde A^t),\, DG(\widetilde A^t)^+\right), \quad \text{where}\quad 
DG(\widetilde A^t):= \varinjlim \widetilde A^t: \Z^{\n}\hookleftarrow
\]
with $DG(\widetilde A^t)^+$ determined (dynamicallly) as before: 

\medskip

Since   $\lambda_\psi$ is simple,  $\widetilde A^t$ has a codimension one invariant subspace $W$ (in $\R^{\n}$) complementary to the eigenspace of $\lambda_\psi$. Let
 \[
 \begin{array}{r}
 (\Z^{\n})^+:=
\{x\in\Z^{\n}\,:\  x\ \text{is in the half-space determined by  $W$ containing}  \\
\text{ the nonnegative eigenvector corresponding to}\  \lambda_\psi \}.
\end{array}
\]
Then $DG(\widetilde A^t)^+:= \left\{[k,g]\in DG(\widetilde A^t)\,:\ g\in (\Z^{\n})^+\right\}$.

We know from \cite{bargesmith} that $\left(DG(\widetilde A^t),\ DG(\widetilde A^t)^+\right)$ is isomorphic to $\left(H^1(\tT_\psi),\ H^1(\tT_\psi)^+\right)$ via an isomorphism that conjugates the shifts.  Our remaining task, then, is to embed $\left(PDG(\psi),PDG(\psi)^+\right)$ into $\left( DG(\widetilde A^t),\ DG(\widetilde A^t)^+\right)$.

Besides enabling us to describe its augmented cohomology as a dimension group, another advantage  of a prepared substitution (like $\psi$)  is that its return lattice $\Sigma$ is the same as 
$\Gamma=\mathrm{pr}_V(\Z^n)$, making formulas for geometric realization simpler. To show these lattices agree, let $\gamma\in \T_\psi$, and let $J$ be an edge of type $j$ in $\gamma$ followed by 
the edge $I$ of type $i$. Since $\psi$ is strictly proper, $\Psi(J)$ begins with an edge of the same type, say $r$, as does $\Psi(I)$. Thus the vector
\[
\begin{array}{r}
\min \Psi(I) - \min\Psi(J)=\max \Psi(J)-\min\Psi(J) =\mathrm{pr}_V(Ae_j)\\ =(A|_V)\mathrm{pr}_Ve_j=(A|_V)\varv_j\in \Theta(j).
\end{array}
\]

 So, $\varv_j\in \Sigma$ for all $j$, and $\Sigma=\Gamma$.  

We may suppose, without loss of generality, that $\{\varv_1,\dots,\varv_d\}$ is a basis for $\Sigma=\Gamma$ , 
taken from  the basis $\{\varv_1,\dots,\varv_n\}$ for $V$. Let $\overline Br(\psi)=\pi_0\varg_\psi(\mathcal C_\psi)$ be the 
branch locus in $V/\Sigma$    (so  $F_{L_\psi}\left(\overline Br(\psi)\right)=Br(\psi)$).  If $\gamma$ and 
$\gamma'$ are asymptotic in $\T_\psi$, then  $\pi_0\varg_\psi(\gamma)$ and $ \pi_0\varg_\psi(\gamma')$                                                                          
are asymptotic in $V/\Sigma$ under the Kronecker flow; that is, $\pi_0\varg_\psi(\gamma)=\pi_0\varg_\psi(\gamma')$. It follows that 
\[
\overline Br(\psi)=\left\{\pi_0\varg_\psi(\gamma_1^f),\dots,\pi_0\varg_\psi(\gamma_{n_f}^f)\right\}\bigcup \left\{
\pi_0\varg_\psi(\gamma_1^b),\dots,\pi_0\varg_\psi(\gamma_{n_b}^b)\right\}.
\]
If $I$ is the unique edge of $\gamma_i^b$ that meets $\E^s$, let $x_i^b:=\min I$; and if $J$ is the unique edge of $\gamma_j^f$ that meets $\E^s$, let $x_j^f:=\max J$.  

Since $\Sigma=\Gamma$, we may take the $u_i=0$ in the definition of geometric realization $\varg_\psi$. Thus $\pi_0\varg_\psi(\gamma_j^f)=x^f_j+\Sigma=:\overline x_j^f$ and $\pi_0\varg_\psi(\gamma_i^b)=x^b_i+\Sigma=:\overline x_i^b$.

If $\#(\overline Br(\psi))=1$,  then 
\[
\left(PDG(\psi),\ PDG(\psi)^+\right)\cong \left(\varinjlim A|_V,\  (\varinjlim A|_V)^+\right),
\]
which is easily embedded in $\left(DG(\psi),\ (DG(\psi)^+\right)$.

Thus we can assume $\#Br(\psi)=m_f+m_b$ with $m_f>0$ and $m_b>0$, and after reindexing, 
\[
\overline Br(\psi)=\left\{\overline x_i^b\,:\ i=1,\dots, m_b\right\}\,\bigcup\, \left\{\overline x_j^f\,:\ j=1,\dots, m_f\right\}.
\]

Recall that $\sigma_i:=\{t\varv_i\,:\ t\in [0,1]\}$ is the oriented segment. The homology classes $[\overline \sigma_i]$ of the oriented cycles $\overline \sigma_i:=\sigma_i+\Sigma,\ i=1,\dots, d$, constitute a basis for $H_1(V/\Sigma)$.

For each $i=1,\dots, n_b$, let $\alpha_i$ denote the oriented line segment in $V$ from $x_i^b$ 
to $0$; for each $j=1,\dots, n_f$, let $\beta_j$ be the directed line segment from $0$ to $x_j^f$; and, let 
$\overline \alpha_i:=\alpha_i+\Sigma,\ \overline \beta_j:=\beta_j+\Sigma$. Then 
\[
\left\{[\os_1],\dots, [\os_\alpha], [\oa_1\ob_1],\, [\oa_1\ob_2],\dots,[\oa_1\ob_{m_f}],\,[\oa_2\ob_1],\,[\oa_3\ob_1],\dots, [\oa_{m_b}\ob_1]\right\}
\]
 is a basis for $\ds H_1\left(V/\Sigma,\,\overline Br(\psi)\right)$.

Let $\L:\Z^\n\to \Z^\n$ (recall $\n=n+n_f+n_b-1$) denote the homomorphism represented by $\widetilde A$ in the standard basis $\mathcal B=\{e_1,\dots, e_\n\}$. For convenience, relabel the basis elements as follows:
\[
e_{ij}=\begin{cases} e_i & \text{for $j=0$ and $i=1,\dots,n$} \\
                            e_{n+j} & \text{for $i=1$ and $j=1,\dots, n_f$} \\
                            e_{n+n_f+i-1} & \text{for $j=2$ and $i=2,\dots,n_b$.}
\end{cases}
\]
We define a homomorphism $\ds P:\Z^\n\to H_1\left(V/\Sigma,\overline Br(\psi)\right)$ on the basis $\mathcal B$ as follows:
\[
\begin{array}{l}
P(e_{i0}):=\sum_{j=1}^d r_{ij}[\os_j],\ \text{provided}\ \varv_i=\sum_{j=1}^d r_{ij}\varv_j,\ i=1,\dots n\text{;}\\[.1in]
P(e_{ij}):=[\oa_i\ob_j],\ i=1,\ j=2,\dots, n_f\ \text{and}\ j=1,\ i=2,\dots, n_b.
 \end{array}
\]
Clearly, $P$ is surjective. The relation  $P\L=(F_{A|_V})_*P$ evidently holds on the basis elements 
$e_{10},\dots, e_{n0}$, so consider $e_{ij}\in\mathcal B$ with $j>0$. Let $p_i:=\llbracket p^b_i\rrbracket\in\Z^n,\ s_j=\llbracket s_j^f\rrbracket\in\Z^n$, and let $\imath:\Z^n\to\Z^\n$ be given by $\imath(x_1,\dots,x_n)=\sum_{i=1}^nx_i\,e_{i0}$. Then, 
$\L(e_{ij})=\imath(p_i)+\imath(s_j)+e_{ij}$ and 
\[
P\L(e_{ij})=P\L(p_i)+P\L(s_j)+P(e_{ij})=\sum_{k=1}^d a_k[\os_k]+\sum_{k=1}^d b_k[\os_k]+[\oa_i\ob)j],
\]
 where $\mathrm{pr}_V(p_i)=\sum_{k=1}^d a_k\varv_k$ and $\mathrm{pr}_V(s_j)=\sum_{k=1}^d b_k\varv_k$.

On the other hand, $(F_{A|_V})_*P(e_{ij})=(F_{A|_V})_*([\oa_i\ob_j])=[\overline{A(\alpha_i\beta_j)}]$, if 
$A(\alpha_i\beta_j)$ denotes the image of the directed  curve $\alpha_i\beta_j$ under the linear map $A:V\to V$. 
We claim that $\overline{A(\alpha_i\beta_j)}$ is homologous to $\overline \rho_i+\oa_i\ob_j+\overline \eta_j$, where $\rho_i:=\{t\mathrm{pr}_V(p_i)\,:\ 0\le t\le 1\}$ and $\eta_j:=\{t\mathrm{pr}_V(s_j)\,:\ 0\le t\le 1\}$ are directed segments. 
\newline
Indeed, this follows
from $A(-x_i^b)=\mathrm{pr}_V(p_i)+(-x^b_i)$ and $Ax_j^f=x_j^f+\mathrm{pr}_V(s_j)$ in the cover $V$ of $V/\Sigma$ (See Figure~\ref{figure1}). Clearly $\overline \rho_i$ is homologous to $\sum_{k=1}^d a_k\os_k$ and $\overline \eta_j$ is homologous to 
$\sigma_{k=1}^d b_k\os_k$. Thus $P\L(e_{ij})=(F_{A|_V})_*P(e_{ij})$.

Since $P$ is surjective, the dual 
\[
P^t:{\rm Hom}\left(H_1\left(V/\Sigma,\overline Br(\psi)\right),\Z\right)\cong H^1\left(V/\Sigma,\overline Br(\psi)\right)\to {\rm Hom}(\Z^\n,\Z)\cong \Z^\n
\]
is injective.  As $\L^tP^t=P^t(F_{A|_V})^*$, we know $P^t$ induces an injection 
\[
\widehat P^t:\varinjlim F_{A|V}^*=PDG(\psi)\to \varprojlim \L^t=DG(\psi) 
\]
that commutes with the shifts. Since the positive cones are determined dynamically, $\widehat P^t$ either  
preserves or reverses order. As $[\sigma_1]^*\in H^1(V/\Sigma,\, \overline Br(\psi))^+$ and $\widehat P^t[\sigma_1]^*=[e_1]^*\in (DG(\psi))^+,\ \widehat P^t$ preserves order. 

Finally, 
\[
\widehat F_{L_\psi}^*:\varinjlim \left( F_{M_\psi}^*:H^1(\mathbb T^d,Br(\psi))\hookleftarrow\right)=PDG(\psi)\to \varinjlim \left(F_{A|_V}^*:H^1\left(V/\Sigma,\overline Br(\psi)\right)\hookleftarrow\right)
\]
by $\widehat F_{L_\psi}^*([(k,g)])=\left[\left(k,\, F_{L_\psi}^*(g)\right)\right]$ is an ordered isomorphism (that commutes with the shifts), which embeds  $\left(PDG(\psi),PDG(\psi)^+\right)$ into $\left(DG(\psi),DG(\psi)^+\right)$.

\small
\begin{figure}
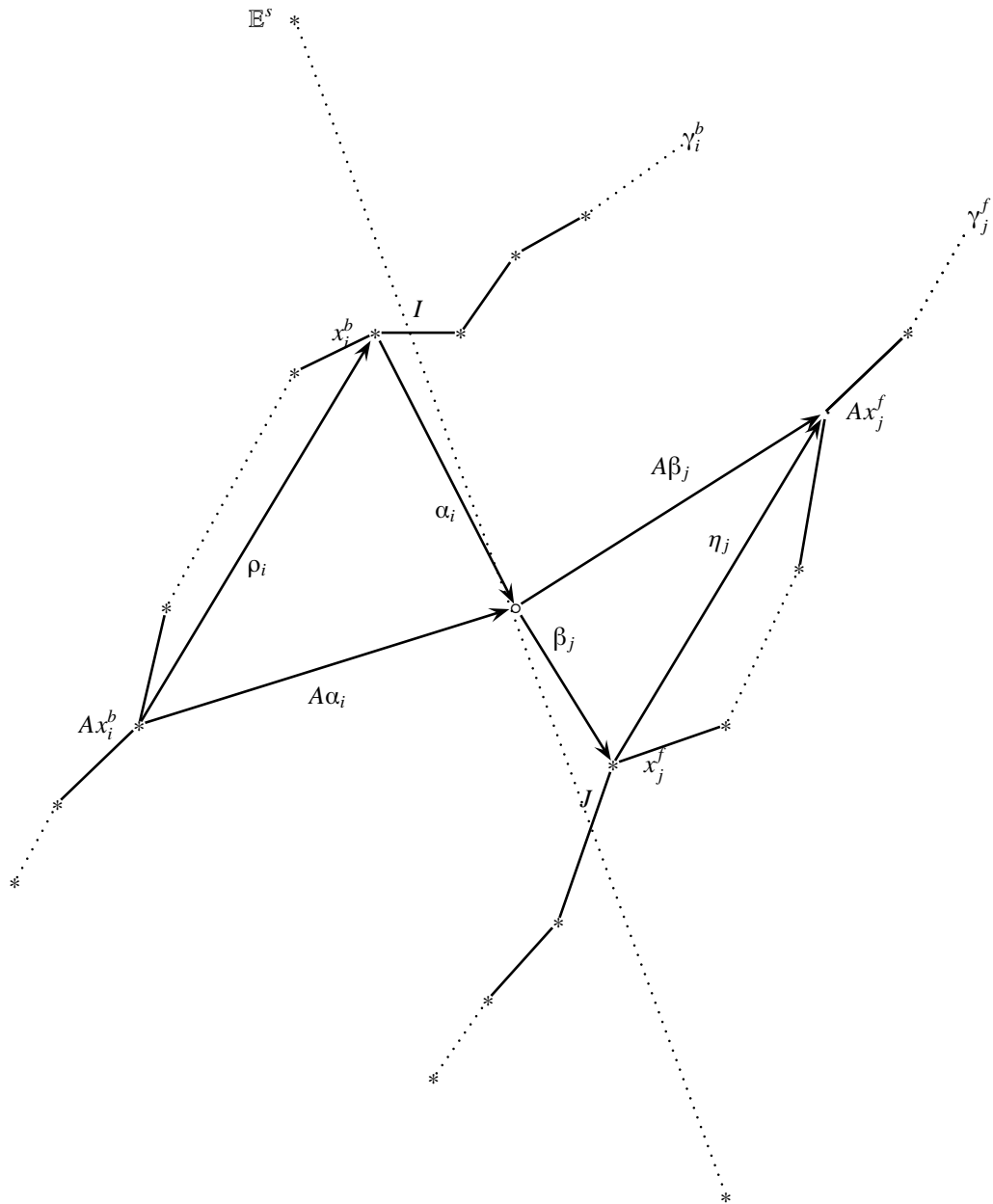
\label{figure1}\caption{A Homology}
\[
  \begin{psmatrix}[mnode=R,rowsep=1pt, colsep=6pt]
\tenands \large{\mathbb E^s} & \small{*}\tenands\tenands\fiveands &\\ 
\allands    \\ 
\allands \\ 
\tenands\tenands\fiveands & \large{\gamma_i^b} \tenands & \\ 
\allands \\    
\tenands\tenands &&& \small{*} \tenands &&&& \large{\ \gamma_j^f}\\ 
\tenands \tenands  \small{*} \tenands \fiveands &&  \\ 
\allands \\ 
\tenands &&& x_i^b  & \small{*} &&&& \small{*} &  \tenands \fiveands  \small{*} &&&\\ 
\tenands & \small{*} \tenands\tenands \fiveands & \\ 
\tenands\tenands\tenands && \,\cdot & \large{Ax_j^f} &&&& \\ 
\allands\\ 
\allands\\ 
\allands \\ 
\tenands\tenands\tenands & \small{*}\fiveands \\ 
\fiveands \small{*} \tenands\fiveands \circ \tenands\fiveands &&\\ 
\allands \\ 
&&&&  \tenands &&& \tenands  \tenands \\ 
&&&\large{Ax_i^b}&\small{*} &\tenands\tenands && \small{*} \tenands    \\ 
\tenands\tenands &&&& \small{*} & \large{\,\,x_j^f}\tenands && \\ 
&& \small{*} && \tenands\tenands &\tenands && \\ 
\allands \\ 
\small{*} \allands \\ 
\tenands\tenands && \small{*} \tenands\fiveands \\ 
\allands\\ 
\tenands\fiveands  &&&& \small{*} \fiveands &&   \tenands  \\ 
\allands\\ 
\tenands \fiveands && \small{*}  \tenands\tenands \\
\allands \\ 
\allands \\
\tenands\tenands \fiveands & & \small{*}  \fiveands &&&&& 
\psset{linestyle=solid, arrowsize=2mm, linewidth=1pt}
\ncline[linestyle=dotted]{6, 24}{4, 27}
\ncline{7, 21}{6, 24}
\ncline{->}{16,21}{20,25}^{\beta_j}
\ncline{->}{20,25}{11,33}^{\eta_j}
\ncline{->}{16,21}{11,33}^{A\beta_j}
\ncline{11,33}{9,35}
\ncline[linestyle=dotted]{9,35}{6,38}
\ncline{->}{19,5}{16,21}_{A\alpha_i}
\ncline{->}{9,15}{16,21}_{\alpha_i}
\ncline{19,25}{18,28}
\ncline[linestyle=dotted]{19,28}{15,32}
\ncline{15,32}{11,33}
\ncline{11,33}{9,35}
\ncline[linestyle=dotted]{9,35}{6,38}
\ncline{9,19}{7,21}
\ncline{10,12}{9,15}
\ncline{9,15}{9,19}^I 
\ncline{->}{19,5}{9,15}_{\rho_i}
\ncline[linestyle=dotted]{1,12}{31,28}
\ncline[linestyle=dotted]{28,18}{26,20}
\ncline{26,20}{24,23}
\ncline{21,3}{19,5}
\ncline{19,5}{16,6}
\ncline[linestyle=dotted]{16,6}{10,12}
\ncline[linestyle=dotted]{23,1}{21,3}
\ncline{24,23}{20,25}^J
\ncline{20,25}{19,28}
\end{psmatrix}
\]
\end{figure}
\normalsize

\end{proof}

{\Small {\parindent=0pt Department of Mathematics, Montana State University,
Bozeman, MT 59717

barge@math.montana.edu, swanson@math.montana.edu

Department of Mathematics, College of Charleston, Charleston, SC 29424
 
diamondb@cofc.edu}}

\end{document}